\documentclass[11pt,reqno]{amsart}

\usepackage{amsmath}
\usepackage{amsxtra}
\usepackage{amscd}
\usepackage{amsthm}
\usepackage{amsfonts}
\usepackage{amssymb}
\usepackage{eucal}
\usepackage[all]{xy}

\numberwithin{equation}{section}

\newtheorem{theorem}{Theorem}[subsection]
\newtheorem{proposition}[theorem]{Proposition}
\newtheorem{lemma}[theorem]{Lemma}
\newtheorem{corollary}[theorem]{Corollary}
\theoremstyle{definition}
\newtheorem{definition}[theorem]{Definition}
\newtheorem{remark}[theorem]{Remark}

\newcommand{\nc}{\newcommand}
\nc{\C}{{\mathbb C}}

\let\al\alpha

\newcommand{\Z}{{\mathbb Z}}

\newcommand{\Ref}[1]{{$($\ref{#1}$)$}}
\newcommand{\bean}{\begin{eqnarray*}}
\newcommand{\eean}{\end{eqnarray*}}
\newcommand{\be}{$$}
\newcommand{\bea}{\begin{eqnarray}}
\newcommand{\ena}{\end{eqnarray}}
\newcommand{\ee}{$$}
\def\ov[#1,#2]{\overset{\scriptstyle #1}{#2}}
\nc{\verm}{M_{k,l}}
\def\wsl{\widehat{\mathfrak{sl}}_2}
\nc{\slt}{{\widehat{\mathfrak{sl}}_2}}
\nc{\inte}{{L_{k,l}}}
\nc{\pkln}{{\mathcal P}_{k,l}^{(N)}}
\nc{\ckln}{{c_{k,l}^{(N)}}}
\def\qbin[#1;#2]{{\left[\matrix{\displaystyle #1}\\{\displaystyle #2}\endmatrix\right]}}
\def\geq{\ge}
\def\leq{\le}
\def\cal{\mathcal}

\nc{\Vk}{{\mathfrak V}_k}

\newcommand\no{\nonumber}

\def\xx[#1,#2]{{#1}^{(#2)}}
\def\yy[#1,#2,#3]{{#1}^{(#2)}_{#3}}

\def\v(#1;#2){\Bigl({#1\atop#2}\Bigr)}

\nc{\alb}{{\boldsymbol{\alpha}}}
\nc{\beb}{{\boldsymbol{\beta}}}

\def\hookdownarrow%
{\Big\downarrow\kern-.14cm\smash{\raise.4cm\hbox{$\scriptstyle\cap$}}}
\def \slg {\mathfrak{sl}_2}

\newcommand{\ch}{{\rm ch\,}}
\nc{\heis}{\widetilde{\mathfrak H}}
\nc{\dualW}{{W_k^{(M,N)}[l_1,l_2,l_3]}^*}
\nc{\dualWa}{{W_k[l_1,l_2,l_3]}^*}
\nc{\card}{{\#}\,}

\begin{document}
\title[Combinatorics of Coinvariants]
{Combinatorics of the $\widehat{sl}_2$ Spaces of Coinvariants III}
\author{B. Feigin, R. Kedem,
 S. Loktev, T. Miwa
 and E. Mukhin}
\address{BF: Landau institute for Theoretical Physics, Chernogolovka,
142432, Russia}
\address{RK: Dept. of Mathematics, 
U. Illinois, Urbana/Champaign}
\address{SL: Institute for Theoretical and Experemental Physics and 
Independent University of Moscow}
\address{TM: Dept, of Mathematics, Kyoto University, Kyoto 606 Japan}
\address{EM: Dept. of Mathematics, Indiana University - Purdue
  University - Indianapolis.}

\begin{abstract}
We give the fermionic character formulas for the spaces
of coinvariants obtained {}from level $k$ integrable representations
of $\widehat{\mathfrak sl}_2$. We establish the functional realization of
the spaces dual to the coinvariant spaces. We parameterize functions
in the dual spaces by rigged partitions, and prove the recursion
relations for the sets of rigged partitions.
\end{abstract}
\maketitle
\section{Introduction}
\subsection{Coinvariant spaces of $\slt$}
Let ${\mathfrak a}$ 
be a Lie subalgebra of a Lie algebra $\mathfrak g$, and $L$ a
$\mathfrak g$-module. The quotient space $L/{\mathfrak a}L$
is called the space of coinvariants of $L$ with respect to ${\mathfrak a}$.
In \cite{FKLMM1,FKLMM2} we studied spaces of
coinvariants for integrable $\widehat{\mathfrak sl}_2$-modules. The present
paper is Part III of the series. We make extensive use of the results
of the pervious papers.

In this paper, we consider the following special case of the
coinvariant. Let $e_i$, $f_i$, $h_i$ ($i\in{\Z}$) be the loop
generators of $\widehat{\mathfrak sl}_2$, and ${\mathfrak
a}={\mathfrak a}^{(M,N)}$ the subalgebra generated by $\{e_i(i\geq
M);f_i(i\geq N)\}$.  Let $L=L_{k,l}$ be the level-$k$ integrable
highest weight $\widehat{\mathfrak sl}_2$-module with highest weight $
(k-l)\Lambda_0+l\Lambda_1$. We are interested in the coinvariant 
\begin{equation}\label{sltchar}
L_{k,l}^{(M,N)}=L_{k,l}/{\mathfrak a}^{(M,N)}L_{k,l}.
\end{equation}
The main result of \cite{FKLMM2} was a theorem about the dimension of
this space, which we showed is given by the Verlinde rule:
\begin{theorem}\label{DIM} For $M,N\geq 0$,         
\begin{equation} 
{\rm dim}\,L_{k,l}^{(M,N)}=\card\Bigl({\cal P}^{M+N}_{k,l}\Bigr),
\end{equation}
where ${\cal P}^N_{k,l}$ is the set of level-$k$ admissible paths of
length $N$ and weight $l$.
(see \cite{FKLMM2} for the precise definition).
\end{theorem}

In fact, the coinvariant space inherits a graded structure from the
integrable module $L_{k,l}$. Let $d$ denote the homogeneous degree
element of $\slt$,
$[d,x_i]=i x_i$ for $x\in{\mathfrak sl}_2$ and define the Hilbert
polynomial or character of the coinvariant space to be
\[
\chi^{(M,N)}_{k,l}(z,q)={\rm trace}_{L^{(M,N)}_{k,l}}q^dz^{h_0}
\]
where $h_0=h\in{\mathfrak sl}_2$. In \cite{FKLMM2} we used a recursion
relation for such characters to prove Theorem \ref{DIM}.  The purpose
of this paper is to derive explicit formulas for these polynomials. It
turns out that our procedure naturally results in fermionic formulas
for the characters. See \cite{FS,St} for some related formulas in the special case $l=0$.

\subsection{The Heisenberg loop algebra and coinvariants}

In order to study the dimension of the coinvariant, in \cite{FKLMM2}
we introduced the simpler coinvariants associated with modules of the
Heisenberg loop algebra.

Let $\mathfrak H$ be the three dimensional Heisenberg algebra with
generators $e,f,h$ and relations $[e,f]=h$ and $h$ central (note that
we use the same notation for the generators of $\slg$, but the
relations are different; it should be clear from the context which
algebra the generators belong to). Let
$\heis$ be the algebra of loops into $\mathfrak H$, generated by
$\{e_i, f_i, h_i; i\in\Z\}$ with relations
\be
[e_i,f_j]=h_{i+j},\qquad [h_i,e_j]=[h_i, f_j]=0.
\ee
Note, that in contrast to $\widehat{\mathfrak sl}_2$,
$\heis$ has a triple-grading, with degrees defined by
\begin{equation}\label{GRADING}
{\rm deg}\,e_i=(1,0,i),\quad{\rm deg}\,f_i=(0,1,i),\quad{\rm deg}\,h_i=(1,1,i).
\end{equation}

Let $W_{k}[l_1,l_2,l_3]$ 
be the $k$-restricted $\heis$-module 
(see \Ref{heis int} for the definition). It is the analog of the
level-$k$ $\slt$-modules, although it is not irreducible. It turns out
that there is a simple relationship between the characters of these
modules and those of $L_{k,l}$.

We consider the coinvariants of $W_k[l_1,l_2,l_3]$ with respect to the
$\heis$ subalgebras ${\mathfrak a}={\mathfrak a}^{(M,N)}$ generated by
the set of elements $\{e_i(i\geq M);f_i(i\geq N)\}$: \be
W_k^{(M,N)}[l_1,l_2,l_3]=W_k[l_1,l_2,l_3]/{\mathfrak
a}^{(M,N)}W_k[l_1,l_2,l_3].  \ee (In this section, we assume
$M,N\geq1$, but in the main text, we treat $M,N\geq 0$.) The
$\heis$-modules and coinvariants inherit the triple-grading
\Ref{GRADING}, and we define the character by \be
\chi^{(M,N)}_k[l_1,l_2,l_3](z_1,z_2,q) =\sum_{m,n,d}{\rm
dim}(W_k^{(M,N)}[l_1,l_2,l_3]_{m,n,d})\; z_1^mz_2^nq^d, \ee where
$W_k^{(M,N)}[l_1,l_2,l_3]_{m,n,d}$ is the subspace of degree
$(m,n,d)$.

In \cite{FKLMM2}, we showed that 
$\widehat{\mathfrak{sl}}_2$-coinvariants and $\heis$-coinvariants
are closely related, and that
$\chi^{(M,N)}_{k,l}$ is given in terms of
$\chi_k^{(M,N)}[l_1,l_2]=\chi_k^{(M,N)}[l_1,l_2,\min(l_1,l_2)]$:
\be
z^l\chi_{k,l}^{(M,N)}(q,z)=
\chi_k^{(M+1,N)}[l,k-l](q^{-2}z^2,z^{-2},q)
-q \chi_k^{(M+1,N)}[l-1,k-l-1](q^{-2}z^2,z^{-2},q).
\ee

The key property of $\chi_k^{(M,N)}[l_1,l_2,l_3]$ used in
the proof of Theorem \ref{DIM}, is that it satisfies the following recursion
relation with respect to $(M,N)$ (see Theorem 6.1.5 of \cite{FKLMM2}:
\begin{theorem}\label{RECREL}
\begin{equation}
\chi^{(M,N)}_k[l_1,l_2,l_3](z_1,z_2,q)
=\sum_{{0\leq a\leq l_3}\atop {0\leq c\leq l_2-a}}
z_1^az_2^{a+c}q^{a+c}\;\chi^{(M,N-1)}_k[l_1',l_2',l_3'](z_1,qz_2,q)
\end{equation}
where $l_1'={\rm min}(l_1+c-a,k-a),\quad l_2'=k-c,\qquad l_3'=l_1'+l_2'-k$.
\end{theorem}
In this paper, we give an explicit formula for the characters
$\chi^{(M,N)}_k[l_1,l_2,l_3]$ (see Theorem \ref{FERMIONIC}).
These formulas have a fermionic form in the sense of 
\cite{KKMM}.

\subsection{Functional realization of dual spaces}
The basic idea in deriving closed forms for the characters is
to consider the function spaces $\dualW$ dual to
$W_k^{(M,N)}[l_1,l_2,l_3]$. The defining relations for $\heis$ are simpler
than those for $\widehat{\mathfrak sl}_2$ because they
respect the grading (\ref{GRADING}). As a consequence, for each fixed
$m,n$, the space $\dualW_{m,n}$ can be realized as a subspace of the
space of rational functions $F(x_1,\dots,x_m;y_1,\dots,y_n)$,
symmetric in each set $\{x_1,\dots,x_m\}$ and $\{y_1,\dots,y_n\}$
separately, having at most simple poles when $x_i=y_j$ and
zeros on the submanifolds $x_i=x_j=y_l$ ($i\not=j$) and $x_i=y_j=y_l$
($j\not=l$).  The dual space is characterized by the vanishing of
functions $F$ on certain submanifolds depending on
$k,l_1,l_2,l_3$. For example, the restriction related to the level $k$
reads as
\begin{equation}\label{VAN}
F=0\quad\text{if $x_1=\dots=x_{k+1}$ or $y_1=\dots=y_{k+1}$}.
\end{equation}
(See section \ref{dual to coinv} for the full definition.)
Because of the high codimensionality of these submanifolds,
it is not possible to immediately deduce the formulas
for the characters, and it is necessary
to introduce a filtration on the dual space,
such that adjoint graded spaces are isomorphic simply to spaces of symmetric
functions, and thus have simple characters. We follow \cite{FS} in this
process.

Let
\begin{equation}\label{MUNU}
\mu=\bigl(k^{m_k},(k-1)^{m_{k-1}},\dots,1^{m_1}\bigr),\quad
\nu=\bigl(k^{n_k},(k-1)^{n_{k-1}},\dots,1^{n_1}\bigr)
\end{equation}
be level-$k$ restricted partitions of $m$ and $n$,
respectively, so that $\sum_\alpha\alpha m_\alpha=m$ and
$\sum_\alpha\alpha n_\alpha=n$.  We consider the following family of
submanifolds
\begin{equation}\label{SUBM}
{\cal M}_{\mu,\nu}:x^{(\alpha)}_{i,1}=\dots=x^{(\alpha)}_{i,\alpha}\quad
(1\leq\alpha\leq k;1\leq i\leq m_\alpha),
\quad y^{(\alpha)}_{i,1}=\dots=y^{(\alpha)}_{i,\alpha}\quad
(1\leq\alpha\leq k;1\leq i\leq n_\alpha),
\end{equation}
where the sets of variables $\{x_j\}$, $\{y_j\}$ are relabeled
$\{x^{(\alpha)}_{i,l}\}$ and
$\{y^{(\alpha)}_{i,l}\}$, respectively.

A subspace ${\cal F}_{\mu,\nu}\subset \dualWa_{m,n}$ is the
 subspace of functions vanishing on 
the submanifolds ${\cal M}_{\mu,\nu}$. Using lexicographic ordering on
partitions, these give a filtration of the dual space, and
the adjoint graded space to this filtration has a
simple structure. For example, if $l_3={\rm min}(l_1,l_2)$,
the graded component corresponding
to $(\mu,\nu)$ is spanned 
by the set of all symmetric polynomials on ${\cal M}_{\mu,\nu}$.
More precisely, we identify the $(\mu,\nu)$-graded
component with the space of functions of the form
$G_{\mu,\nu}g$,
where $G_{\mu,\nu}$ is a fixed rational function depending only on
 $\mu,\nu$ and $g$ is an arbitrary
polynomial in the variables
$\{x^{(\alpha)}_i\}_{1\leq i\leq m_\alpha}$ and
$\{y^{(\alpha)}_i\}_{1\leq i\leq n_\alpha}$, $1\leq\alpha\leq k$,
 symmetric under the exchange of variables with the same superscript $\alpha$, 
$x^{(\alpha)}_i\leftrightarrow x^{(\alpha)}_j$ or 
$y^{(\alpha)}_i\leftrightarrow y^{(\alpha)}_j$.
This space has a basis Sym($\prod_{\alpha,i} 
(x_i^{(\alpha)})^{r_i^{(\alpha)}}(y_i^{(\alpha)})^{s_i^{(\alpha)}}$),
where for each $\alpha$,
$r^{(\alpha)}=\{r^{(\alpha)}_{1\leq i\leq m_\alpha}\}$
and
$s^{(\alpha)}=\{s^{(\alpha)}_{1\leq i\leq n_\alpha}\}$ are
sets of integers satisfying
$r^{(\alpha)}_1\geq\cdots\geq r^{(\alpha)}_{m_\alpha}\geq0$ and
$s^{(\alpha)}_1\geq\cdots\geq s^{(\alpha)}_{n_\alpha}\geq0$, respectively.

These basis elements are in one to one correspondence with
combinatorial data $(\mu,r;\nu,s)$ called {\it rigged partitions},
introduced in \cite{KKR,KR}.
The set of non-negative integers $r$ is called a rigging of the partition
$\mu$.

If $l_3<\min(l_1,l_2)$, there is an additional restrictions for the riggings 
from below,
\begin{equation}\label{TAURES}
r^{(\alpha)}_i+s^{(\beta)}_j\geq \min(\alpha,\beta)-\max(\al-l_1,0)-
\max(\beta-l_2,0)-l_3.
\end{equation}

The space dual to the coinvariant, $\dualW$, is the subspace of functions $F$
 which satisfy the degree restrictions
\[
{\rm deg}_{x_1}F<M,\quad{\rm deg}_{y_1}F<N.
\]
We will show that the degree restrictions
translates to conditions for the riggings $r$
and $s$ of the form
\begin{equation}\label{DEGRES}
r^{(\alpha)}_i\leq P^{(M)}_{\mu,\nu}[l_1]_\al,\quad
s^{(\alpha)}_i\leq Q^{(N)}_{\mu,\nu}[l_2]_\al,
\end{equation}
where the vacancy numbers
$P^{(M)}_{\mu,\nu}[l_1]$, $ Q^{(N)}_{\mu,\nu}[l_2]$ are defined in
equations \Ref{p vac}, \Ref{q vac}.

Our final result is that the adjoint graded space of 
$\dualW_{m,n}$ has a basis
labeled by pairs of rigged partitions $(\mu,r;\nu,s)$ with the
restrictions on the riggings of the form (\ref{TAURES}) and
(\ref{DEGRES}). Denote the set of such rigged partitions by
$R^{(M,N)}_{m,n}[l_1,l_2,l_3]$.  Because of Theorem \ref{RECREL}, one
can expect that there is an inductive construction of
$R^{(M,N)}_{m,n}[l_1,l_2,l_3]$ {}from
$R^{(M,N-1)}_{m-a,n-a-c}[l_1',l_2',l_1'+l_2'-k]$.  In fact, this is
true, and we will describe it explicitly.

The logical ordering of this paper is somewhat different.
We prove directly that the evaluation map which maps the space of
functions of the form $F$ to the space of functions spanned by $G_{\mu,\nu}g$
is injective. However, we do not have a simple direct proof that it is
surjective. Instead, we construct $R^{(M,N)}_{m,n}[l_1,l_2,l_3]$
inductively, and this assures the surjectivity by dimension counting arguments.

The plan of paper is as follows. In Section \ref{RIG} we give preliminaries on
rigged partitions and state the main recursion theorem (Theorem
\ref{recursion theorem}). 
In Section \ref{FUN}, we construct the functional realization of dual spaces,
their filtrations and describe the adjoint graded spaces.
We also give the resulting fermionic formulas for the characters.
Sections \ref{RIG PREL}, \ref{lower summation section} and
\ref{upper summation section} are devoted to the proof of Theorem
\ref{recursion theorem}.
The arguments in these sections are purely combinatorial.
In Section \ref{RIG PREL},
we define admissible pairs $(I,J)$ of index sets
belonging to $\{1,\dots,k\}$. Then, we define two types of subsets
of rigged partitions indexed by admissible pairs, the lower and upper subsets.
We construct a bijection {}from the upper to the lower subsets indexed by
the same pair $(I,J)$. In Sections \ref{lower summation section} and
\ref{upper summation section}, we give the decompositions
of the set of rigged partitions for $(M,N)$ by the lower subsets, and
that for $(M,N-1)$ by the upper subsets, respectively. This 
completes the proof of Theorem \ref{recursion theorem}.

\bigskip

\noindent
{\it Acknowledgments.}\quad 
The work of SL is in part supported by grant RFBR-01-01-00546

\section{Rigged partitions and the main recursion theorem}\label{RIG}
We define level restricted rigged partitions, and
state the main recursion theorem for sets of rigged partitions, Theorem
\ref{recursion theorem}, together with an outline the proof. 

\subsection{Rigged partitions and vacancy numbers}
Let $k\in\Z_{\geq 1}$, $m\in\Z_{\geq0}$ and $I_k=\{1,2,...,k\}$. Let $\mu$ be
a level-$k$ restricted partition of $m$, that is
\begin{equation}
\mu=(k^{m_k},\dots,2^{m_2},1^{m_1}),  \qquad \sum_{\al=1}^k\al m_\al=m,
\label{YD}
\end{equation}
We denote by $m_\alpha(\mu)$ the number of rows of length $\alpha$ in
the partition (or Young diagram) $\mu$.
%

A rigging of $\mu$ is a set of integers 
$r=\{r_i^{(\al)}\}_{\al\in I_k,\ 1\leq i\leq m_\al(\mu)}$ such that
\begin{equation}
r_1^{(\al)}\geq \dots \geq r_{m_\al(\mu)}^{(\al)}\geq0\quad(\al\in I_k).
\label{R1}
\end{equation}
A partition with a rigging, $(\mu,r)$, is called a rigged partition.
Denote by $R_m$ the set of all such level-$k$ restricted rigged
partitions of $m$. We set $R_{m,n}=R_m\times R_n$. 


Let $l_1,l_2,l_3$ be integers satisfying
\begin{equation}
0\leq l_1,l_2\leq k,\quad0\leq l_3\leq {\rm min}(l_1,l_2).
\label{l condition}
\end{equation}
Define
\begin{eqnarray}
\tau^{(\al,\beta)}[l_1,l_2,l_3]&=&
\min(\al,\beta,l_1,l_2,l_1+\beta-\al,
l_2+\al-\beta,l_1+l_2-\al,l_1+l_2-\beta)-l_3\no\\
&=&
{\rm min}(\alpha,\beta)-(\alpha-l_1)^+-(\beta-l_2)^+-l_3,\label{NEWTAU}
\end{eqnarray}
where 
$x^+={\rm max}(x,0),\quad x^-={\rm max}(-x,0).$
We define a subset of $R_{m,n}$ where the lower bounds of the riggings
are restricted by \Ref{NEWTAU}:
\begin{eqnarray}
R_{m,n}[l_1,l_2,l_3]=\Bigl\{(\mu,r;\nu,s)\in R_{m,n}\ ;\ 
r_{m_\al(\mu)}^{(\al)}+s_{m_\beta(\mu)}^{(\beta)}\geq\tau^{(\al,\beta)}[l_1,l_2,l_3]\Bigr\}
\label{exception condition}
\end{eqnarray}
Since $\tau^{(\alpha,\beta)}[l_1,l_2,{\rm min}(l_1,l_2)]\leq0$,
there is no restriction in this case, and
$R_{m,n}[l_1,l_2,{\rm min}(l_1,l_2)]=R_{m,n}$

Let $M,N$ be non-negative integers. Define vectors of vacancy numbers 
$P_{\mu,\nu}^{(M)}[l_1], Q_{\mu,\nu}^{(N)}[l_2]\in\Z^k$, where
\begin{align}
&P_{\mu,\nu}^{(M)}[l]_\al=
\al M-(\al-l)^+
+\sum_{\beta=1}^k\min(\al,\beta)(m_\beta(\nu)-2m_\beta(\mu)),\label{p vac}\\
&Q_{\mu,\nu}^{(N)}[l]_\al= P_{\nu,\mu}^{(N)}[l]_\al
.\label{q vac}
\end{align}

We define the subset $R_{m,n}^{(M,N)}[l_1,l_2]\subset R_{m,n}$:
\begin{eqnarray}
R_{m,n}^{(M,N)}[l_1,l_2]=\{(\mu,r;\nu,s)\in R_{m,n}&:&
P_{\mu,\nu}^{(M)}[l_1]_\al,Q_{\mu,\nu}^{(N)}[l_2]_\al\geq0,
\label{COC}\\
&&r^{(\alpha)}_1\leq P_{\mu,\nu}^{(M)}[l_1]_\al,\
s^{(\alpha)}_1\leq Q_{\mu,\nu}^{(N)}[l_2]_\al
\label{M degree condition}
\}.
\end{eqnarray}
The first condition, (\ref{COC}), is
non-trivial only in the case $m_\al(\mu)=0$ or $m_\al(\nu)=0$. Otherwise,
it follows from (\ref{M degree condition}).
However, see Proposition \ref{CUTPRO} for the actual implication
of this conditions.

Finally define the set
\begin{equation}
R_{m,n}^{(M,N)}[l_1,l_2,l_3]=R_{m,n}^{(M,N)}[l_1,l_2]\cap R_{m,n}[l_1,l_2,l_3].
\end{equation}
It is defined
for negative values of $m,n$ by
\be
R_{m,n}^{(M,N)}[l_1,l_2,l_3]=\emptyset,\qquad m<0\;\; {\rm or}\;\; n<0.
\ee

Before passing, we prove
\begin{proposition}\label{CUTPRO}
Suppose that the conditions $(\ref{M degree condition})$ 
hold. Then,
the conditions $(\ref{COC})$ are equivalent to the following requirements:
\begin{eqnarray}
\text{If $M=0$ then $n-2m\geq k-l_1$};\label{M=0}\\
\text{If $N=0$ then $m-2n\geq k-l_2$.}\label{N=0}
\end{eqnarray}
Namely, it is enough to require the conditions $(\ref{COC})$ only for the cases
$M=0,\alpha=k$ and $N=0,\alpha=k$.
\end{proposition}
\begin{proof}
In the following, when we write a condition concerning the $0$-th
component of a $k$ vector
(e.g., the case $i=0$ for $P_{\mu,\nu}^{(M)}[l_1]_i\geq0$
in the next paragraph or $P_i\geq\rho_i$ in the proof
of Lemma \ref{REFCOM}), we mean that the condition is void.

First we prove that if $M\geq1$
the condition  $P_{\mu,\nu}^{(M)}[l_1]_k\geq0$
follows {}from (\ref{M degree condition}).
Suppose otherwise, there exists
$0\leq i\leq k-1$ such that $P_{\mu,\nu}^{(M)}[l_1]_i\geq0$
and $m_{i+1}=\dots=m_k=0$.  (Here,
$m_\al=m_\al(\mu),n_\al=m_\al(\nu)$.)
Then, we have
\begin{eqnarray*}
0&>&P_{\mu,\nu}^{(M)}[l_1]_k\notag\\
&=&P_{\mu,\nu}^{(M)}[l_1]_i+(k-i)M+(i-l_1)^+-(k-l_1)
+\sum_{\beta\geq i+1}(\beta-i)n_\beta\notag\\
&\geq&0,
\end{eqnarray*}
which is a contradiction. Similarly, if $N\geq1$ the condition 
$Q_{\mu,\nu}^{(N)}[l_2]_k\geq0$ follows {}from (\ref{M degree condition}).

Now we will prove that for all $M\geq0$,
the conditions $P_{\mu,\nu}^{(M)}[l_1]_\al\geq0$
for $1\leq\alpha\leq k-1$ follow {}from $P_{\mu,\nu}^{(M)}[l_1]_k\geq0$.
(The proof is similar for $Q_{\mu,\nu}^{(N)}[l_2]_\al$.)

Suppose otherwise, there exists $i$ and $j$ such that
$0\leq i<j-1\leq k-1$, $P_{\mu,\nu}^{(M)}[l_1]_i\geq0$,
$P_{\mu,\nu}^{(M)}[l_1]_j\geq0$ and $P_{\mu,\nu}^{(M)}[l_1]_\alpha<0$
(and thereby $m_\alpha=0$) for $i+1\leq\alpha\leq j-1$.
Set $p=\frac1{j-i}$ so that we have $pi+(1-p)j=j-1$.
Then we have
\begin{eqnarray}
0&\leq& pP_{\mu,\nu}^{(M)}[l_1]_i+(1-p)P_{\mu,\nu}^{(M)}[l_1]_j\notag\\
&=&P_{\mu,\nu}^{(M)}[l_1]_{j-1}+(j-1-l_1)^+-p(i-l_1)^+-(1-p)(j-l_1)^+\notag\\
&&\quad +
\sum_\beta\Bigl(p\;{\rm min}(i,\beta)+(1-p){\rm min}(j,\beta)
-{\rm min}(j-1,\beta)\Bigr)n_\beta
\notag\\
&=&P_{\mu,\nu}^{(M)}[l_1]_{j-1}-\theta(i+1\leq l_1\leq j-1)(l_1-i)p
-\sum_{i+1\leq\beta\leq j-1}(\beta-i)pn_\beta\notag\\
&<&0,\label{REFCOM2}
\end{eqnarray}
which is a contradiction. Here we used the notation
\begin{eqnarray}
\theta(*)=
\begin{cases}
1\quad\text{ if $*$ is true;}\\
0\quad\text{ if $*$ is false.}\\
\end{cases}
\end{eqnarray}
\end{proof}

Proposition \ref{CUTPRO} implies
\begin{corollary}
For $(M,N)=(0,0)$ we have
\begin{equation}\label{INIT}
R^{(0,0)}_{m,n}[l_1,l_2,l_3]=
\begin{cases}
\{(\emptyset,\emptyset;\emptyset,\emptyset)\}
\quad&\text{if $l_1=l_2=k$ and $m=n=0$;}\\
\emptyset\quad&\text{otherwise}.
\end{cases}
\end{equation}
\end{corollary}
\subsection{Recursion Theorem for rigged partitions}

We state the main theorem on recursion.

\begin{theorem}\label{recursion theorem} The cardinalities of the sets of the rigged partitions satisfy the following relation$:$
\be
\card(R_{m,n}^{(M,N)}[l_1,l_2,l_3])=\sum_{0\leq a\leq l_3 \atop
0\leq c\leq l_2-a}\card
(R_{m-a,n-a-c}^{(M,N-1)}[l_1',l_2',l_3']),
\ee
where
\bea
l_1'&=&l_1+c-a-(l_1+c-k)^+,\notag\\
l_2'&=&k-c,\label{primes}\\
l_3'&=&l_1'+l_2'-k.\notag
\ena
\end{theorem}
In what follows, we fix the notation $l_1',l_2',l_3'$
to be the integers given by (\ref{primes}), and $b=a+c$.
Theorem \ref{recursion theorem} is proved in Sections \ref{RIG PREL},
\ref{lower summation section} and \ref{upper summation section}.

Let us outline the idea of the proof.
We construct an explicit bijection 
\be
{\mathfrak m}:\;\bigsqcup_{{0\leq a\leq l_3}\atop{0\leq c\leq l_2-a} }
R_{m-a,n-a-c}^{(M,N-1)}[l_1',l_2',l_3']\to R_{m,n}^{(M,N)}[l_1,l_2,l_3]
\ee
in several steps.
In Section \ref{LOWSEC}, for $I,J\subset\{1,\dots,k\}$, we define the subsets
$R_{m,n}^{(M,N)}[l_1,l_2]_{I,J}\subset R_{m,n}$ (the lower subsets).
In Section \ref{UPSEC}, we define
$R_{m-a,n-b}^{(M,N-1)}[l_1]^{I,J}\subset R_{m-a,n-b}$ (the upper subsets),
where $a=\card(I)$ and $b=\card(J)$. In Section \ref{BIJECTION},
we construct the bijection
\be
{\mathfrak m}_{I,J}:\;R_{m-a,n-b}^{(M,N-1)}[l_1]^{I,J}\to
R_{m,n}^{(M,N)}[l_1,l_2]_{I,J}.
\ee
In Section \ref{lower summation section} we will prove that 
for each $(l_1,l_2,l_3)$ satisfying (\ref{l condition})
\be
R_{m,n}^{(M,N)}[l_1,l_2,l_3]=\bigsqcup_{{I,J\subset\{1,\dots,k\}}
\atop{\card(I)\leq l_3,\;\;\card(J)\leq l_2}}
R_{m,n}^{(M,N)}[l_1,l_2]_{I,J},
\ee
and in Section \ref{upper summation section} that
for each $(l_1,a,c)$ and $(l'_1,l'_2,l'_3)$ determined by (\ref{primes})
\be
R_{m-a,n-b}^{(M,N-1)}[l'_1,l'_2,l'_3]=\bigsqcup_{{I,J\subset\{1,\dots,k\}}
\atop{\card (I)=a, \card (J)=b}
}R_{m-a,n-b}^{(M,N-1)}[l_1]^{I,J}.
\ee
This will complete the proof of Theorem \ref{recursion theorem}.

An important implication of Theorem \ref{recursion theorem}, and the
main interest we have in proving it, is the following result.
\begin{corollary}
Fix an integer $k\in{\Z}_{\geq1}$, and consider the spaces of coinvariants
of the $\heis$-modules, $W_k^{(M,N)}[l_1,l_2,l_3]$ (see
\Ref{heis int},\Ref{heis inte})
and the sets of rigged partitions
\be
R^{(M,N)}[l_1,l_2,l_3]=\sqcup_{m,n}R^{(M,N)}_{m,n}[l_1,l_2,l_3].
\ee
Then
\begin{equation}\label{DIMSH}
{\rm dim}\,W_k^{(M,N)}[l_1,l_2,l_3]=\card(R^{(M,N)}[l_1,l_2,l_3]).
\end{equation}
\end{corollary}
\begin{proof}
Using Theorem \ref{RECREL} with $z_1=z_2=q=1$, we see that
these two sets of numbers satisfy the same recursion with the same
initial condition.
\end{proof}

\section{Functional realization of dual spaces and character formulas}\label{FUN}
In this section we identify the space dual to the module
$W_k[l_1,l_2,l_3]_{m,n}$
with a certain space of rational functions in $m+n$ variables.
We introduce a filtration in this space and describe the adjoint graded space
explicitly by using the rigged partitions. As a corollary we compute the
character of the space of coinvariants $W^{(M,N)}_k[l_1,l_2,l_3]$.

\subsection{Dual of the universal enveloping algebra}

Let $\heis$ be the Heisenberg loop algebra with
generators $e_i,f_i,h_i$ $(i\in \Z$) and relations
\be
[e_i,f_j]=h_{i+j},\qquad [e_i,h_j]=[f_i,h_j]=0.
\ee
Consider its universal enveloping algebra $U\heis$. The algebra
$U\heis$ is graded by
\[
{\rm deg}\,e_i=(1,0),\quad{\rm deg}\,f_i=(0,1),\quad{\rm deg}\,h_i=(1,1).
\]
Let $(U\heis)_{m,n}$ be the subspace of degree $(m,n)$.
We construct the space dual to $(U\heis)_{m,n}$
in the space of rational functions in the variables
$(x_1,\ldots,x_m;y_1,\ldots,y_n)$.

Consider the space of rational functions
\begin{eqnarray}
{\cal F}_{m,n}&=&\{F=\frac p{\prod_{i,j}(x_i-y_j)}\ :\ 
p\in \C[x_1^{\pm1},\ldots,x_m^{\pm1},y_1^{\pm1},\ldots,y_n^{\pm1}],
\nonumber\\
&&\text{symmetric in $x_1,\ldots,x_m$ and $y_1,\ldots,y_n$ separately,}
\nonumber\\
&&\text{where $p=0$ if $x_1=x_2=y_1$ or $x_1=y_1=y_2$.}\}\label{P}
\end{eqnarray}

There exists a coupling between $(U\heis)_{m,n}$ and 
${\cal F}_{m,n}$. In order to define it, consider
the mappings $L_{e_i}:{\cal F}_{m,n}\rightarrow{\cal F}_{m-1,n}$,
$L_{f_i}:{\cal F}_{m,n}\rightarrow{\cal F}_{m,n-1}$,
$L_{h_i}:{\cal F}_{m,n}\rightarrow{\cal F}_{m-1,n-1}$:
\begin{eqnarray}
L_{e_i}(F)&=&\oint\frac{dx_1}{2\pi \sqrt{-1}x_1}Fx_1^{-i},\label{L1}\\
L_{f_i}(F)&=&\oint\frac{dy_1}{2\pi \sqrt{-1}y_1}Fy_1^{-i},\label{L2}\\
L_{h_i}(F)&=&\oint\frac{dy_1}{2\pi \sqrt{-1}y_1}\Bigl\{(x_1-y_1)F\Bigr\}
\Bigl|_{x_1=y_1}y_1^{-1-i},\label{L3}
\end{eqnarray}
where $F\in {\cal F}_{m,n}$.
In each of these equations, we take the contour of integration to be a circle
in ${\C}$ oriented counter-clockwise such that all the poles are inside.
Because of the vanishing of $p$ at $x_1=x_2=y_1$ and $x_1=y_1=y_2$,
the integrand of (\ref{L3}) has the only pole in $y_1$ at $y_1=0$.

Similarly, we define the mappings $R_{e_i}$, $R_{f_i}$, $R_{h_i}$
by the same formulas (\ref{L1}), (\ref{L2}), (\ref{L3}), respectively,
using a contour such that all the poles except the origin are outside.
As we noted above, we have $L_{h_i}=R_{h_i}$.


The following proposition is standard. We omit the proof.
\begin{proposition}\label{CPL}
There exists a unique coupling between
$(U\heis)_{m,n}$ and ${\cal F}_{m,n}$ such that
\begin{eqnarray*}
&&\langle1,1\rangle=1\\
&&\langle e_iw,F\rangle=\langle w,L_{e_i}(F)\rangle,\quad
\langle f_iw,F\rangle=\langle w,L_{f_i}(F)\rangle,\quad
\langle h_iw,F\rangle=\langle w,L_{h_i}(F)\rangle,\\
&&\langle we_i,F\rangle=\langle w,R_{e_i}(F)\rangle,\quad
\langle wf_i,F\rangle=\langle w,R_{f_i}(F)\rangle,\quad
\langle wh_i,F\rangle=\langle w,R_{h_i}(F)\rangle.
\end{eqnarray*}
\end{proposition}
For example, it follows immediately that
\begin{lemma}\label{easy coupling}
If $w=e_{i_1}\dots e_{i_m}f_{j_1}\dots f_{j_n}$, then the coupling 
$\langle w,F\rangle$ is equal to the coefficient of
$x_1^{i_1}\dots x_m^{i_m}y_1^{j_1}\dots y_n^{j_n}$ in the Laurent series
obtained by expanding $F$ in positive powers of $y_j/x_i$.
\end{lemma}

\begin{proposition}\label{nondegenerate}
The coupling given by Proposition \ref{CPL} is non-degenerate.
\end{proposition}
\begin{proof}
First we show that for any nonzero
$F\in{\cal F}_{m,n}$ there exists $w\in (U\heis)_{m,n}$
such that $\langle w,F\rangle\not=0$.

Consider the lexicographic ordering of monomials
$x_1^{i_1}\dots x_m^{i_m}y_1^{j_1}\dots y_n^{j_n}$. Namely,
the monomial $x_1^{i_1}\dots x_m^{i_m}y_1^{j_1}\dots y_n^{j_n}$
is higher than $x_1^{i'_1}\dots x_m^{i'_m}y_1^{j'_1}\dots y_n^{j'_n}$
if $i_1>i'_1$, or if $i_1=i'_1$ and $i_2>i'_2$, and so on. Let
$x_1^{i_1}\dots x_m^{i_m}y_1^{j_1}\dots y_n^{j_n}$ be the highest
monomial present in $p$ of $F$ in (\ref{P}). Then, taking
$w=e_{i_1-n}\dots e_{i_m-n}f_{j_1}\dots f_{j_n}$ we have
$\langle w,F\rangle\not=0$.

Next we show that for any nonzero
$w\in (U\heis)_{m,n}$ there exists $F\in{\cal F}_{m,n}$
such that $\langle w,F\rangle\not=0$.
For $l\leq {\rm min}(m,n)$ let $Z_l$ be the set of indices
$({\bf k},{\bf i},{\bf j})$ such that
${\bf k}\in\Z^{l}$, ${\bf i}\in\Z^{m-l}$ and ${\bf j}\in\Z^{n-l}$,
with
$$
k_1\leq\dots\leq k_l,\ i_1\leq\dots\leq i_{m-l},\ 
j_1\leq\dots\leq j_{n-l}.
$$
By the PBW theorem the monomials
\begin{equation}\label{MONO}
M[{\bf k,i,j}]:=
h_{k_1}\dots h_{k_l}e_{i_1}\dots e_{i_{m-l}}f_{j_1}\dots f_{j_{n-l}}
\end{equation}
span $(U\heis)_{m,n}$.
Set
\[
F[{\bf k,i,j}] 
={\rm Sym}\Bigl(\prod_{a=1}^l\frac{y_a^{k_a}}{x_a-y_a}
\prod_{b=1}^{m-l}x_{l+b}^{i_b}\prod_{c=1}^{n-l}y_{l+c}^{j_c}\Bigr).
\]
Take
$({\bf k},{\bf i},{\bf j})\in Z_l$ and
$({\bf k'},{\bf i'},{\bf j'})\in Z_{l'}$.
Using the definition of the coupling, we have
\begin{eqnarray}
\langle
M[{\bf k,i,j}],F[{\bf k',i',j'}]\rangle =
\begin{cases}
0&\text{if $l>l'$};\\
\delta_{\bf k,k'}\delta_{\bf i,i'}\delta_{\bf j,j'}&\text{if $l=l'$}.
\end{cases}
\end{eqnarray}
The assertion follows {}from this.
\end{proof}

\subsection{Dual to the $\heis$-module and coinvariant}\label{dual to coinv}

Following \cite{FKLMM2}, define the $\heis$-module
$W[l_1,l_2,l_3]$ as a quotient of $U\heis$
by the left ideal generated by the elements
\begin{equation}\label{heis int}
x_i\,\,(i\leq0, x\in {\mathfrak H}),\
e_1^{l_1+1},\quad f_1^{l_2+1},\quad h_1^{l_3+1},\
\end{equation}
and the level-$k$ restricted module $W_k[l_1,l_2,l_3]$
is the quotient of $W[l_1,l_2,l_3]$ by the two-sided ideal generated by
\begin{equation}\label{heis inte}
e(z)^{k+1},\quad f(z)^{k+1},
\end{equation}
where we used the generating series
$e(z)=\sum_{i\in\Z} e_iz^i$, $f(z)=\sum_{i\in\Z}f_iz^i$. (Strictly
speaking, these elements are in the completion of $U\heis$; however as
usual, the module is in the category $\cal O$ due to \Ref{heis int}
and, when
acting in $W[l_1,l_2,l_3]$, they are finite sums in $U\heis$.)

The dual space of $W_k[l_1,l_2,l_3]_{m,n}$ is realized in 
${\cal F}_{m,n}$ as the subspace orthogonal to
these ideals. We denote this subspace by
$\dualWa_{m,n}$. The following theorem is 
a consequence of Proposition \ref{CPL}
\begin{theorem}\label{dualspace}
The space $\dualWa_{m,n}$ is given by
\begin{eqnarray*}
\dualWa_{m,n}
&=&\{F=\frac{\prod_ix_i\prod_jy_j}{\prod_{i,j}(x_i-y_j)}f\in{\cal F}_{m,n};\\
&&\text{$f\in \C[x_1,...,x_m,y_1,...,y_n]$ },\\
&&\text{$f=0$ if $x_1=\dots=x_{k+1}$ or $y_1=\dots=y_{k+1}$or}\\
&&\text{$x_1=\dots=x_{l_1+1}=0$ or $y_1=\dots=y_{l_2+1}=0$, }\\
&&\text{$\prod_{i=1}^{l_3}\left(\frac{\partial}{\partial x_{i+1}}\frac{\partial}{\partial y_{i+1}}\right)^if=0 $ if}\\
&&x_1=\dots=x_{l_3+1}=y_1=\dots=y_{l_3+1}=0.\}
\end{eqnarray*}
\end{theorem}
Take $M,N\in{\Z}_{\geq0}$. Let the subalgebra
${\mathfrak a}^{(M,N)}$ of $\heis$
be generated by the elements $e_i$ ($i\geq M$) and $f_i$ ($i\geq N$).
Following \cite{FKLMM2}, define the space of coinvariants by
\begin{eqnarray}
W^{(M,N)}_k[l_1,l_2,l_3]=\oplus_{m,n}
W^{(M,N)}_k[l_1,l_2,l_3]_{m,n},\label{COI}\\
W^{(M,N)}_k[l_1,l_2,l_3]_{m,n}=
\begin{cases}
0\quad\text{if $M=0$ and $n-2m<k-l_1$;}\\
0\quad\text{if $N=0$ and $m-2n<k-l_2$;}\\
W_k[l_1,l_2,l_3]/{\mathfrak a}^{(M,N)}W_k[l_1,l_2,l_3],\quad\text{otherwise}.
\end{cases}\nonumber
\end{eqnarray}

Define ${\cal F}_{m,n}^{(M,N)}\subset{\cal F}_{m,n}$ to be the subset
consisiting of functions $F$ satisfying the degree restrictions
\begin{equation}\label{DEG}
{\rm deg}_{x_1}F<M,\quad{\rm deg}_{y_1}F<N.
\end{equation}
Here the degree of the rational function $F$ in the variable $x_1$ is
defined to be the highest power in $x_1$ appearing in the Laurent series
expansion of $F$ in positive powers of $y_j/x_1$. In other words,
we have ${\rm deg}_{x_1}F=1-n+{\rm deg}_{x_1}f$. Similarly, we have
${\rm deg}_{y_1}F=1-m+{\rm deg}_{y_1}f$. If $m$ or $n$ is zero, the
corresponding degree restriction is void.

\begin{definition}{\rm
We define the space of rational functions $\dualW$
by
\begin{eqnarray}
\dualW&=&\oplus_{m,n}\dualW_{m,n},
\nonumber\\
\dualW_{m,n}&=&
\begin{cases}
0\quad\text{if $M=0$ and $n-2m<k-l_1$;}\\
0\quad\text{if $N=0$ and $m-2n<k-l_2$;}\\
W_k[l_1,l_2,l_3]^*_{m,n}\cap{\cal F}^{(M,N)}_{m.n}\quad\text{otherwise}.
\end{cases}
\end{eqnarray}
}
\end{definition}
The space $\dualW_{m,n}$ is finite-dimensional
and dual to $W^{(M,N)}_k[l_1,l_2,l_3]_{m,n}$ given by (\ref{COI}).

\subsection{Polynomials with Serre relations}
In this section we study symmetric polynomials of the form
$f(x_1,\dots,x_m;y_1,\dots,y_n)$ 
which vanish when $x_1=x_2=y_1$ or $x_1=y_1=y_2$. Proposition
\ref{h1prop} will be used in the proof of Theorem \ref{FILISO}
(see Lemma \ref{ZL3}).

For a function $f(x_1,\dots,x_m;y_1,\dots,y_n)$ and 
$a_i,b_j\in{\Z}_{\geq 0}$ denote
\be
[a_1,\dots,a_s;b_1,\dots,b_t]f:=\prod_{i=1}^s\Bigl(\frac{\partial}{\partial
  x_i}\Bigr)^{a_i} 
\prod_{i=1}^t\Bigl( \frac{\partial}{\partial
  y_i}\Bigr)^{b_i}f|_{x_1=\dots=x_s=y_1=\dots=y_t=z}.
\ee
We also denote the functions
\begin{align*}
&\Bigl(
    \frac{\partial}{\partial z}\Bigr)^r
\prod_{i=3}^s\Bigl (\frac{\partial}{\partial x_i}\Bigr)^{a_i}
\prod_{i=2}^t\Bigl (\frac{\partial}{\partial
    y_i}\Bigr)^{b_i}f(z,z,x_3,\dots;z,y_2,\dots) 
    |_{x_3=\dots=x_s=y_2=\dots=y_t=z}, \\
& \Bigl(
    \frac{\partial}{\partial z}\Bigr)^r
\prod_{i=2}^s\Bigl(\frac{\partial}{\partial x_i}\Bigr)^{a_i}
\prod_{i=3}^t\Bigl(\frac{\partial}{\partial y_i}\Bigr)^{b_i}f(z,x_2,\dots;z,z,y_3,\dots)
    |_{x_2=\dots=x_s=y_3=\dots=y_t=z}
\end{align*}
by
\be
[a_3,\dots,a_s,(*,*;*)^r,b_2,\dots,b_t]f,\qquad 
[a_2,\dots,a_s,(*;*,*)^r,b_3,\dots,b_t]f.
\ee

If the number of $x$ variables in $f$ is smaller than the $s$ or the
number of $y$ variables is smaller than $t$ then we define the functions
$[a_1,\dots,a_s;b_1,\dots,b_t]f$,
$[a_3,\dots,a_s,(*,*;*)^r,b_2,\dots,b_t]f$,

\noindent $[a_2,\dots,a_s,(*;*,*)^r,b_3,\dots,b_t]f$  to be $0$.  

We have relations
\begin{align}
&[a_3,\dots,a_s,(*,*;*)^r,b_2,\dots,b_t]f=\sum_{a_1+a_2+b_1=r}
\frac{r!}{a_1!a_2!b_1!}[a_1,\dots,a_s;b_1,\dots,b_t]f,\label{tr1}\\
&[a_2,\dots,a_s,(*;*,*)^r,b_3,\dots,b_t]f=\sum_{a_1+b_1+b_2=r}
\frac{r!}{a_1!b_1!b_2!}[a_1,\dots,a_s;b_1,\dots,b_t]f.\label{tr2}
\end{align}

For the rest of this section, 
let $f(x,y)$ be a polynomial in $x_1,\dots,x_m,y_1,\dots,y_n$,
symmetric with respect to permutations of $x$ and to permutations of
$y$, satisfying the {\it Serre} relations: 
\bea\label{Serre}
f(x,y)=0 \qquad {\rm if}\;\; x_1=x_2=y_1 \quad {\rm or}\;\; x_1=y_1=y_2. 
\ena
Note that now $[a_1,\dots,a_s;b_1,\dots,b_t]f$ does not depend on the
order of $a_i$ or $b_j$. Also
we have
\be
[a_3,\dots,a_s,(*,*;*)^r,b_2,\dots,b_t]f=[a_2,\dots,a_s,(*;*,*)^r,b_3,
\dots,b_t]f=0.  
\ee
In particular we 
have many linear relations among $[a_1,\dots,a_s;b_1,\dots,b_t]f$
thanks to \Ref{tr1}, \Ref{tr2}. The following lemma describes some
of the relations which consist of a single term. 

\begin{lemma} For $s,t\in\Z_{\geq 0}$, we have the identities
\begin{align}
&[0,1,2,\dots,s,s,s^t;0,1,2,\dots,s-1,s+t]f=0,\label{serre1}\\
&[0,1,2,\dots,s-1,s+t;0,1,2,\dots,s,s,s^t]f=0.\label{serre2}
\end{align}
where $s^t$ denotes $s,s\dots,s$ repeated $t$ times.
\end{lemma}
\begin{proof}

We use induction on $s$. 

The identity \Ref{serre1} for $s=0$ takes the form
\bea\label{serre trivial}
[0^{t+2};t]f=0,\qquad t\in\Z_{\geq 0}, \qquad 
\ena
The case $t=0$ is just the Serre relation: $[0,0;0]f=0$. 
We obtain \Ref{serre trivial}
by induction on $t$. Assume $[0^{t+2};t]f=0$,
$t=0,\dots,t_0-1$. Then $[0^{t_0+2};t_0]f=0$ follows from 
the identity $[0^{t_0}(*,*;*)^{t_0}]f=0$. Indeed on the RHS of
\Ref{tr1} for $[0^{t_0}(*,*;*)^{t_0}]$,
the only term left is exactly $[0^{t_0+2};s_0]f=0$.  
The $s=0$ case of identity \Ref{serre2} is proved similarly.

Now assume \Ref{serre1}, \Ref{serre1} are 
proved for $s=0,\dots,s_0-1$ and let us prove them
for $s=s_0$. It is enough to prove \Ref{serre1}, then \Ref{serre2} 
is done by the same argument switching the roles of $x$ and
$y$. 

We use induction on $t$. 
The case $t=0$ follows from the
identity 
\be
[0,1,2,\dots,s_0-1,(*,*;*)^{3s_0},0,1,2,\dots,s_0-1]f=0.
\ee 
Suppose we have the statement for $t=0,\dots,t_0-1$, then the case $t=t_0$
follows from the identity
\be
[0,\dots,s_0-1,s_0^{t_0},(*,*;*)^{3s_0+t_0},0,1,\dots,s_0-1]f=0.
\ee 
\end{proof}

Now we derive more identities under additional assumptions. 

For a function $g(x_1,\dots,x_m;y_1,\dots,y_n)$ and
$a_i,b_j\in\Z_{\geq 0}$, we denote
\begin{align*}
&[a_1,\dots,a_s;b_1,\dots,b_t]'g=[a_1,\dots,a_s;b_1,\dots,b_t]g|_{z=0}\\
&[a_3,\dots,a_s,(*,*;*)^r,b_2,\dots,b_t]'g=
[a_3,\dots,a_s,(*,*;*)^r,b_2,\dots,b_t]g|_{z=0},\\ 
&[a_2,\dots,a_s,(*;*,*)^r,b_3,\dots,b_t]'g=
[a_2,\dots,a_s,(*;*,*)^r,b_3,\dots,b_t]g|_{z=0}.
\end{align*}
Then the $h_1^{l_3+1}=0$ relation is translated into
\bea\label{h1}
[0,1,\dots,l_3;0,1,\dots,l_3]'f=0, 
\ena
see Theorem \ref{dualspace}.

\begin{remark}
The condition $e_1^{l_1+1}=0$ reads $[0^{l_1+1};\emptyset]'f=0$. It
follows from our results that we automatically have
$h_1^{l_1+1}=0$. It is an instructive exercise to prove  
$[0,1,\dots,l_1;0,1,\dots,l_1]'f=0$ starting from
$[0^{l_1+1};\emptyset]'f=0$ and using \Ref{tr1}, \Ref{tr2}.
\end{remark}

\begin{proposition}\label{h1prop} Let $f$ satisfy \Ref{h1}.
Then for $s\in\Z_{\geq l_3}$ 
we have the identity
\bea\label{h3cor}
[0^{s+1};s^{l_3+1}]'f=0.
\ena
\end{proposition}
\begin{proof}
We use the induction on $s$. Assume the statement is proved for
$s=l_3,\dots,s_0-1$. (We assume nothing if $s_0=l_3$.) We will prove it
for $s=s_0$. To do that we prove by 
the inverse induction on $r$ the identity 
\bea\label{ax}
[0,1,\dots,r-1,(r)^{s_0-r+1};0,1\dots,r-1,(s_0)^{l_3-r+1}]'f=0,
\ena
where $r=l_3+1,l_3,l_3-1,\dots,0$. The case $r=0$ is exactly \Ref{h3cor} for
$s=s_0$.

The identity \Ref{ax} for $r=l_3+1$ follows directly from
\Ref{h1}. Assume we have \Ref{ax} for $r=l_3+1,l_3,\dots,r_0+1$. Let
us prove it for $r=r_0$. For that we prove the identity
\be
[0,1,\dots,r_0,(r_0)^q(r_0+1)^{s_0-r_0-q};0,1\dots,r_0-1,r_0+q,(s_0)^{l_3-r_0}]'f=0,
\ee
for $q=0,\dots,s_0-r_0$ by induction on $q$. For $q=0$ we have exactly
\Ref{ax} for $r=r_0+1$. If the statement is proved for
$q=0,\dots,q_0-1$ then the statement for $q=q_0$ follows from the
relation
\be
[0,1,\dots,r_0-1,(r_0)^{q_0-1}(r_0+1)^{s_0-r_0-q_0}(**;*)^{3r_0+q_0},
0,1\dots,r_0-1,(s_0)^{l_3-r_0}]'f=0
\ee
and \Ref{serre1}.

For $q=s_0-r_0$ we obtain \Ref{ax} for $r=r_0$ and the proof is finished.
\end{proof}

\subsection{Multiplication of functional spaces}
In this section we describe a multiplicative structure which relates
the functional spaces for different levels $k$. Though the results of this
section are not used in what follows, we think that
Theorem \ref{star prod surj}, is interesting in its own right.

Fix $k^{(j)},\;l_i^{(j)}$ $(i=1,2,3,\;j=1,2$) and set
$k=k^{(1)}+k^{(2)}$, $l_i=l_i^{(1)}+l_i^{(2)}$  $(i=1,2,3)$.

Let $\Delta:\; U\heis\to U\heis\otimes U\heis$
be the usual comultiplication defined by the rule $\Delta(g)=1\otimes
g+g\otimes 1$ for $g\in \heis$. We also denote by $\Delta$ the
map of $ U\heis$ modules
\be
\Delta:\; W_k[l_1,l_2,l_3]\to W_{k^{(1)}}[l_1^{(1)},l_2^{(1)},l_3^{(1)}]\otimes W_{k^{(2)}}[l_1^{(2)},l_2^{(2)},l_3^{(2)}]
\ee
uniquely determined by the condition $\Delta(v)=v^{(1)}\otimes
v^{(2)}$, where $v$, $v^{(1)}$ and $v^{(2)}$ are the highest weight
vectors of the corresponding modules.

The map $\Delta$ descends to the spaces of coinvariants
\be
\Delta^{(M,N)}\; W_k^{(M,N)}[l_1,l_2,l_3]\to W^{(M,N)}_{k^{(1)}}[l_1^{(1)},l_2^{(1)},l_3^{(1)}]\otimes W_{k^{(2)}}^{(M,N)}[l_1^{(2)},l_2^{(2)},l_3^{(2)}].
\ee
By Proposition 6.3.3 in \cite{FKLMM2} the map $\Delta^{(M,N)}$ is injective.

Define the map 
\be
*:\;W^{*}_{k^{(1)}}[l_1^{(1)},l_2^{(1)},l_3^{(1)}]\otimes
W_{k^{(2)}}^{*}[l_1^{(2)},l_2^{(2)},l_3^{(2)}]\to W^{*}_k[l_1,l_2,l_3]
\ee 
by the following rule. Let $F^{(j)}(x_1^{(j)},\dots,x_{m^{(j)}}^{(j)};
y_1^{(j)},\dots,y_{n^{(j)}}^{(j)}\in
W^{*}_{k^{(j)}}[l_1^{(j)},l_2^{(j)},l_3^{(j)}]_{m^{(j)},n^{(j)}}$,
($j=1,2$).
Then $F^{(1)}*F^{(2)}\in W^{*}_k[l_1,l_2,l_3]_{m,n}$,
where $m=m^{(1)}+m^{(2)}$, $n=n^{(1)}+n^{(2)}$,
is given by
\bean
\lefteqn{F^{(1)}*F^{(2)}(x_1,\dots, x_m;y_1,\dots, y_n)=} \\
&&
={\rm Sym}\left(F^{(1)}(x_1,\dots,x_{m^{(1)}};y_1,\dots,y_{n^{(1)}})\;
F^{(2)}(x_{m^{(1)}+1},\dots,x_{m};y_{n^{(1)}+1},\dots,y_{n})\right).
\eean
Here ${\rm Sym}$ denotes the symmetrization with respect to two groups
of variables $x_1,\dots,x_m$ and $y_1,\dots,y_n$. 

\begin{lemma}\label{star dual delta} The map $*$ is well defined. Moreover, the map $*$ is
dual to the map $\Delta$:
\bea\label{dual maps}
\langle \Delta (w), F^{(1)}\otimes F^{(2)}\rangle=\langle w,  F^{(1)}
* F^{(2)}\rangle,
\ena
where $w\in W_k[l_1,l_2,l_3]$ and the pairing on the tensor product of vector spaces is standard:
$\langle v^{(1)}\otimes v^{(2)}, F^{(1)}\otimes F^{(2)}\rangle=
\langle v^{(1)}, F^{(1)}\rangle\langle  v^{(2)},F^{(2)}\rangle$.
\end{lemma}
\begin{proof}
The fact that the map $*$ is well defined follows directly from the
defintion. 
Note that the vectors $w$ of the form $w=e_{i_1}\dots
e_{i_m}f_{j_1}\dots f_{j_n}v$, where $v$ is the highest weight vector,
span $W_k[l_1,l_2,l_3]_{m,n}$. Indeed, as shown in the proof of
Proposition \ref{nondegenerate}, the orthogonal complement of the span
of such vectors is trivial. Therefore it is enough to check \Ref{dual
maps} for $w$.  For such vectors the equation \Ref{dual maps} is
clear {}from  Lemma \ref{easy coupling}.
\end{proof}

The map $*$ obviously descends to the spaces dual to the coinvariants:
\be
*^{(M,N)}:\;W^{*(M,N)}_{k^{(1)}}[l_1^{(1)},l_2^{(1)},l_3^{(1)}]\otimes
W_{k^{(2)}}^{*(M,N)}[l_1^{(2)},l_2^{(2)},l_3^{(2)}]\to
W^{*(M,N)}_k[l_1,l_2,l_3].
\ee 

{}From Lemma \ref{star dual delta} and the injectivity of the coproduct,
Proposition 6.3.3 in
\cite{FKLMM2}, we obtain
\begin{theorem}\label{star prod surj}
The map $*^{(M,N)}$ is surjective.
\end{theorem}
This is a rather simple statement for certain spaces of symmetric functions.
However, we do not know of any direct proof of this statement.

\subsection{Filtration of $W^*_k[l_1,l_2,l_3]$}
Let $\mu$ be a level-$k$ restricted partition of $m$ of the form
\Ref{YD}. We will define a map $\varphi_\mu$ which
sends functions of the variables $(x_1,\ldots,x_m)$
to functions of the variables
$\{x_j^{(\al)}\}_{\al\in I_k,1\leq j\leq m_\al(\mu)}$.

Fix a numbering from $1$ to $m$ of the set of indices $(\alpha,j)$ where
$\alpha\in I_k$ and $1\leq j\leq m_\al(\mu)$.
We define $\varphi(x_i)=x^{(\alpha)}_j$ where $(\alpha,j)$
is the $i$-th index in this numbering. The $\mu$-evaluation
map $\varphi_\mu$ is defined by
\[
\varphi_\mu\ \bigl(F(x_1,...,x_m)\bigr)\ =
F(\varphi_\mu(x_1),...,\varphi_\mu(x_m)).
\]
If $F$ is a symmetric function, then $\varphi_\mu(F)$ is symmetric
in the variables
$(x^{(\alpha)}_1,\dots,x^{(\alpha)}_{m_\alpha})$ with fixed $\alpha$.
Moreover, $\varphi_\mu(F)$ is independent of the choice of the numbering.

Given a pair of partitions $(\mu,\nu)$ of $(m,n)$,
$(\mu,\nu)$-evaluation $\varphi_{\mu,\nu}$ 
is defined by
\[
\varphi_{\mu,\nu}(F(x_1,\dots,x_m;y_1,\dots,y_n))
=F(\varphi_\mu(x_1),...,\varphi_\mu(x_m);
\varphi_\nu(y_1),...,\varphi_\nu(y_n)).
\]

Partitions are ordered lexicographically,
$\mu>\mu'$ if and only if there exists some $i$ for which
$\mu_i>\mu_i'$ and $\mu_j=\mu_j'$ for all $j<i$.
Similarly, pairs  of partitions $(\mu,\nu)$ are ordered,
$(\mu,\nu)>(\mu',\nu')$ if and only if $\mu>\mu'$, or $\mu=\mu'$ and
$\nu>\nu'$.

Now suppose $F\in W^*_k[l_1,l_2,l_3]$.
Since $F$ does not have a pole at $x_i=x_j$ or $y_i=y_j$,
the $(\mu,\nu)$-evaluation is well-defined. Consider the subspaces
\begin{eqnarray}
{\rm Ker}\,\varphi_{\mu,\nu}&\subset&
W^*_k[l_1,l_2,l_3],\label{KER}\\
\Gamma_{\mu,\nu}&=&\cap_{(\mu',\nu')>(\mu,\nu)}
{\rm Ker}\,\varphi_{\mu,\nu},\no\\
\Gamma'_{\mu,\nu}&=&\Gamma_{\mu,\nu}\cap{\rm Ker}\,\varphi_{\mu,\nu}.\no
\end{eqnarray}
The subspaces $\Gamma_{\mu,\nu}$ give a filtration of 
$W^*_k[l_1,l_2,l_3]$. Our goal is to characterize the adjoint graded space 
${\rm Gr}_{\mu,\nu}=\Gamma_{\mu,\nu}/\Gamma'_{\mu,\nu}$
(see Theorem \ref{FILISO}).


\begin{lemma}\label{ZERO}
Let $F\in\Gamma_{\mu,\nu}$. The function $\varphi_{\mu,\nu}(F)$ has a zero
of order at least $2\,{\rm min}(\alpha,\beta)$
if $x^{(\alpha)}_i=x^{(\beta)}_j$ or $y^{(\alpha)}_i=y^{(\beta)}_j$.
\end{lemma}
\begin{proof}
Consider the case $x^{(\alpha)}_i=x^{(\beta)}_j$ with
$\alpha\geq\beta$, with $\al,\beta$ fixed. 
Denote the variables $x_k$ such that $\varphi_\mu(x_k)=x^{(\beta)}_j$
by $x^{(\beta)}_{j,l}$ ($l=1,\ldots,m_\beta$) in some ordering.

We can carry out the evaluation in two steps:
$\varphi_{\mu,\nu}(F)=\varphi_2(\varphi_1(F))$, where $\varphi_1$
is the evaluation of all the variables except
$x^{(\beta)}_{j,l}$ ($l=1,\ldots,m_\beta$) 
and $\varphi_2$ is the evaluation of the variables $x^{(\beta)}_{j,l}$.
Let $F_1=\varphi_1(F)$.
Since $\alpha\geq\beta$ and $F\in\Gamma_{\mu,\nu}$, we have
\[
F_1\Bigl|_{x^{(\beta)}_{j,l}=x^{(\alpha)}_i}=0,\qquad 1\leq l \leq \beta.
\]
Differentiating the left hand side of this equality by $x^{(\alpha)}_i$
and using
the symmetry of $F$ with respect to $(x_1,,\dots,x_m)$, we can deduce that
\[
\frac{\partial F_1}{\partial x^{(\beta)}_{j,l}}
\Biggl|_{x^{(\beta)}_{j,l}=x^{(\alpha)}_i}=0,\qquad 1\leq l \leq \beta.
\]
Therefore, $F_1$ has a zero of order at least two at
$x^{(\beta)}_{j,l}=x^{(\alpha)}_i$ for each $l$. After evaluation,
$\varphi_{2}(F_1)$ is divisible by
$(x^{(\alpha)}_i-x^{(\beta)}_j)^{2\beta}$.
\end{proof}

\begin{lemma}\label{POLE}
Let $F\in\Gamma_{\mu,\nu}$. The function $\varphi_{\mu,\nu}(F)$ has a pole
of order at most ${\rm min}(\alpha,\beta)$ if $x^{(\alpha)}_i=y^{(\beta)}_j$.
\end{lemma}
\begin{proof}
Without loss of generality, we can assume that $\alpha\geq\beta$. If
$\alpha=1$ the assertion follows immediately.

Suppose $\alpha\geq2$. Set
$g=\varphi_{\mu,\nu}(f)$, where $f$ is the polynomial function of
Theorem \ref{dualspace}. It is enough to show that $g$ is divisible by
$(x^{(\alpha)}_i-y^{(\beta)}_j)^{(\alpha-1)\beta}$, because the
evaluation of the prefactor in Theorem
\ref{dualspace} only contains a pole of order $\al\beta$ at this point.

Let
$Y_{[j,\beta]}=\{y^{(\beta)}_{j,l}\}_{1\leq
l\leq\beta}=\varphi_{\mu,\nu}^{-1}(y_i^{(\beta)})$. We obtain $g$ in
two steps: $\varphi_{\mu,\nu}(f)=\varphi_2(\varphi_1(f))$, where
$\varphi_1$ is the evaluation of all the variables except those in
$Y_{[j,\beta]}$. 

Using the fact that $f=0$ if $x_1=x_2=y_1$,
$$\frac{\partial^s}{\partial (x^{(\alpha)}_i)^s}\,\varphi_1(f)
\Bigl|_{y^{(\beta)}_{j,l}=x^{(\alpha)}_i}\Bigr. = 0,
\qquad
0\leq s\leq\alpha-2,\quad 1\leq l\leq\beta.
$$
Therefore, $\varphi_1(f)$
 is divisible by $(y^{(\beta)}_{j,l}-x^{(\alpha)}_i)^{\alpha-1}$,
and hence $g$ is divisible by
$(x^{(\alpha)}_i-y^{(\beta)}_j)^{(\alpha-1)\beta}$.
\end{proof}

\begin{lemma}\label{ZL12}
Let $F\in\Gamma_{\mu,\nu}$ and $f$ be as in Theorem
\ref{dualspace}. The function $\varphi_{\mu,\nu}(f)$ has a zero of
order at least $(\alpha-l_1)^+$ (resp., $(\alpha-l_2)^+$) if
$x^{(\alpha)}_i=0$ (resp., $y^{(\alpha)}_i=0$).
\end{lemma}
\begin{proof}
The assertion follows by a similar argument as in the proof of Lemma
\ref{POLE} {}from the restriction on $f$ that it is zero if
$x_1=\dots=x_{l_1+1}=0$ or $y_1=\dots=y_{l_2+1}=0$.
\end{proof}

Let $f(x,y)$ be a polynomial in two variables $x$ and $y$.
We say that $f$ has a zero of order $s$ at $x=y=0$ if $f(tx,ty)$
has a zero of order $s$ at $t=0$.


\begin{lemma}\label{ZL3}
Let $F\in\Gamma_{\mu,\nu}$ and $f$ be as in Theorem \ref{dualspace}.
Then the function $\varphi_{\mu,\nu}(f)$ has a zero
of order at least $\alpha\beta-l_3$ at
$x^{(\alpha)}_i=y^{(\beta)}_j=0$.
\end{lemma}
\begin{proof}
If $l_3\geq {\rm min}(\alpha,\beta)$ then there is nothing to prove
due to Lemma \ref{POLE}. Therefore, without loss of generality we assume
$l_3+1\leq\alpha\leq \beta$. Let
\begin{align*}
h:=f(x_1,\dots,x_\alpha;\underbrace{y,\dots,y}_\beta), \qquad
g:= \left(\prod_{i=1}^\alpha(x_i-y)^{\beta-1}\right)^{-1}h.
\end{align*}
Note  that for $i=1,\dots,\alpha$, we have
\be
\left(\frac{\partial}{\partial y}\right)^sh|_{x_i=y}=0,\qquad
  s=0,\dots,\beta-2,
\ee
because $f=0$ if $x_i=y_j=y_t$. Therefore, $g$ is a polynomial.   

From Proposition \ref{h1prop} and \Ref{serre trivial} we obtain
\bea\label{sbeta}
\Bigl(\frac{\partial}{\partial y}\Bigr)^{s\beta-l_3-1}h|_{x_1=\dots=x_s=y=0}
=0, \qquad s=l_3+1,l_3+2,\dots,\alpha-1.
\ena

Now, it follows by induction on $r$ that for $r=0,\dots,\alpha-l_3$ the
polynomial $g$ is of the form 
\be
g=y^rg'_r +\sum_{i=0}^{r-1}y^ig_i,
\ee
where $g_r'$ is a polynomial and 
$g_i$ are polynomials independent on $y$ of degree at least 
$\alpha-l_3-i$ in $x_1,\dots,x_\alpha$.
Indeed, if we have the statement for $r=r_0-1$, then
the case $r=r_0$ follows from \Ref{sbeta} with $s=r_0+l_3$. 

Therefore $g$ is of degree at least $\alpha-l_3$ in
$x_1,\dots,x_\alpha,y$ and the lemma follows.
\end{proof}

Let $\mu$ and $\nu$ be level $k$ partitions of $m$ and $n$, respectively. Set
\[
G_{\mu,\nu}=
\prod_{\alpha,i}\Bigl(x^{(\alpha)}_i\Bigr)^{\alpha+(\alpha-l_1)^+}
\prod_{\alpha,i}\Bigl(y^{(\alpha)}_i\Bigr)^{\alpha+(\alpha-l_2)^+}
\frac{\prod(x^{(\alpha)}_i-x^{(\beta)}_j)^{2{\rm min}(\alpha,\beta)}
\prod(y^{(\alpha)}_i-y^{(\beta)}_j)^{2{\rm min}(\alpha,\beta)}}
{\prod(x^{(\alpha)}_i-y^{(\beta)}_j)^{{\rm min}(\alpha,\beta)}}.
\]

Consider the space of rational functions in the variables $\{x^{(\alpha)}_i\}$,
$\{y^{(\alpha)}_i\}$ defined as follows:
\begin{eqnarray}
{\cal G}_{\mu,\nu}[l_1,l_2,l_3]&=&\{G=G_{\mu,\nu}g\ ;\ 
g\in\C[\{x^{(\alpha)}_i\},\{y^{(\alpha)}_i\}],
\nonumber\\
&&\text{$g$ is invariant by the transposition
$x_i^{(\al)}\leftrightarrow x_j^{(\al)}$ or
$y_i^{(\al)}\leftrightarrow y_j^{(\al)}$,}\nonumber\\
&&\text{$g$ has a zero at $x^{(\alpha)}_i=y^{(\beta)}_j=0$ of order at least
$\tau^{(\alpha,\beta)}[l_1,l_2,l_3]$ given by (\ref{NEWTAU}).}\}\nonumber\\
\end{eqnarray}
We define the total homogeneous degree of $G$ as the
homogeneous degree of $G$ in all the variables $x^{(\alpha)}_i$
and $y^{(\alpha)}_i$.
\begin{proposition}\label{well-defined}
The evaluation map
\[
\varphi_{\mu,\nu}:
W^*_k[l_1,l_2,l_3]\rightarrow{\cal G}_{\mu,\nu}[l_1,l_2,l_3],\quad
F\mapsto\varphi_{\mu,\nu}(F)
\]
is well-defined,  injective and preserves the total homogeneous degree.
\end{proposition}
\begin{proof}
This follows {}from Lemmas \ref{ZERO}, \ref{POLE}, \ref{ZL12} and \ref{ZL3}.
\end{proof}

\begin{lemma} If $F$ belongs to the subspace
$\dualW_{m,n}\cap\Gamma_{\mu,\nu}$
defined by the conditions (\ref{DEG}) and (\ref{KER})
then the function $g$ given by $\varphi_{\mu,\nu}(F)=G_{\mu,\nu}g$
satisfies the degree restrictions
\begin{equation}
{\rm deg}_{x^{(\alpha)}_i}g\leq P^{(M)}_{\mu,\nu}[l_1]_\alpha,\quad
{\rm deg}_{y^{(\alpha)}_i}g\leq Q^{(N)}_{\mu,\nu}[l_2]_\alpha,
\end{equation}
where $ P^{(M)}_{\mu,\nu}[l_1], Q^{(N)}_{\mu,\nu}[l_2]$ are given by
\Ref{p vac}, \Ref{q vac}.
\end{lemma}
The proof is straightforward.

Set
\begin{eqnarray}
{\cal G}_k^{(M,N)}[l_1,l_2,l_3]
&=&\oplus_{m,n}{\cal G}^{(M,N)}_{m,n}[l_1,l_2,l_3],\nonumber\\
{\cal G}^{(M,N)}_{m,n}[l_1,l_2,l_3]
&=&\{G=G_{\mu,\nu}g\in{\cal G}_{\mu,\nu}[l_1,l_2,l_3];\nonumber\\
&&|\mu|=m,\quad|\nu|=n,\nonumber\\
&&P^{(M)}_{\mu,\nu}[l_1]\geq0,\quad Q^{(N)}_{\mu,\nu}[l_2]\geq0,\nonumber\\
&&{\rm deg}_{x^{(\alpha)}_i}g\leq P^{(M)}_{\mu,\nu}[l_1]_\alpha,\quad
{\rm deg}_{y^{(\alpha)}_i}g\leq Q^{(N)}_{\mu,\nu}[l_2]_\alpha.\}\nonumber\\
\end{eqnarray}

Consider the filtration of $\dualW_{m,n}$
consisting of the subspaces $\dualW_{m,n}\cap\Gamma_{\mu,\nu}$,
and the adjoint graded space ${\rm Gr}(W^{*(M,N)}_k[l_1,l_2,l_3]_{m,n})$.
We set
\[
{\rm Gr}(W^{*(M,N)}_k[l_1,l_2,l_3])
=\oplus_{m,n}{\rm Gr}(W^{*(M,N)}_k[l_1,l_2,l_3]_{m,n})
\]
The mappings $\varphi_{\mu,\nu}$ induce an injective map
\begin{equation}\label{VARPHI}
\varphi:{\rm Gr}(W^{*(M,N)}_k[l_1,l_2,l_3])\rightarrow
{\cal G}_k^{(M,N)}[l_1,l_2,l_3]
\end{equation}

In fact, this map is an isomorphism, however we do not know of a
straightforward proof of the surjectivity. Nevertheless we can prove
the following theorem:
\begin{theorem}\label{FILISO}
The mapping $\varphi$ of (\ref{VARPHI}) is an isomorphism preserving
the total homogeneous degree.
\end{theorem}
\begin{proof}
For a rigged partition $(\mu,r)$ we denote by $m_r(\{x^{(\alpha)}_i\})$
the monomial symmetric polynomial corresponding to the monomial
$\prod_{\alpha,i}\Bigl(x^{(\alpha)}_i\Bigr)^{r^{(\alpha)}_i}$.
The space of rational functions
${\cal G}^{(M,N)}_{m,n}[l_1,l_2,l_3]$
can be parameterized by the set of rigged partitions
$R^{(M,N)}_{m,n}[l_1,l_2,l_3]$ by associating
$G_{\mu,\nu}g$ to $(\mu,r;\nu,s)$
where $g=m_r(\{x^{(\alpha)}_i\})m_s(\{y^{(\alpha)}_i\})$.
The statement follows {}from the injectivity of $\varphi$ and
the equality of dimensions (\ref{DIMSH}).
\end{proof}

\subsection{Characters of coinvariants}
The purpose of this section is to compute
the characters of $W^{(M,N)}_k[l_1,l_2,l_3]$ and the space of
$\slt$-coinvariants. 

The algebra $U\heis$ has a triple grading given by \Ref{GRADING}.
The spaces $W^{(M,N)}_k[l_1,l_2,l_3]$ are quotients of $U\heis$ and
have an induced grading on them. Define the characters of 
$W^{(M,N)}_k[l_1,l_2,l_3]$ to be
\[
\ch_{z_1,z_2,q}W^{(M,N)}_k[l_1,l_2,l_3]
=\sum_{m,n,d}{\rm dim}(W^{(M,N)}_k[l_1,l_2,l_3]_{m,n,d})\;z_1^mz_2^nq^d ,
\]
where $W^{(M,N)}_k[l_1,l_2,l_3]_{m,n,d}$ is the subspace of degree $(m,n,d)$.

Note that the dual space $\dualW$ is similarly graded, with
\[
{\rm deg}\,x_i=(1,0,1),\quad{\rm deg}\,y_i=(0,1,1).
\]
Hence we can define the character of the function spaces described
above. These are equal to those of the corresponding quotients of $U\heis$.

The evaluation mapping preserves the degree. Hence, in the image of
the evaluation by $\varphi$, the induced degree is
\[
{\rm deg}\,x^{(\alpha)}_i=(1,0,1),\quad{\rm deg}\,y^{(\alpha)}_i=(0,1,1).
\]

We can rephrase this in terms of rigged partitions.
Define the degree of a pair of rigged partitions $(\mu,r;\nu,s)$ to be
\[
d(\mu,r;\nu,s)={\rm deg}\,G_{\mu,\nu}+\sum_{\alpha,i} r^{(\alpha)}_i
+\sum_{\alpha,i} s^{(\alpha)}_i,
\]
and the character of the set $R^{(M,N)}_{m,n}[l_1,l_2,l_3]$ by
\begin{equation}\label{RCHAR}
{\ch}_{q}R^{(M,N)}_{m,n}[l_1,l_2,l_3]=\sum_{(\mu,r;\nu,s)
\in R^{(M,N)}_{m,n}[l_1,l_2,l_3]}q^{d(\mu,r;\nu,s)}.
\end{equation}
By definition
\[
\ch_q{\cal G}_{m,n}^{(M,N)}[l_1,l_2,l_3]
=\ch_qR^{(M,N)}_{m,n}[l_1,l_2,l_3].
\]
Finally, by Theorem \ref{FILISO}, we have
\[
\ch_qW^{(M,N)}_k[l_1,l_2,l_3]_{m,n}
=\ch_q{\cal G}_k^{(M,N)}[l_1,l_2,l_3].
\]

Let us compute these characters explicitly. Set
\[
A_{\alpha,\beta}={\rm min}(\alpha,\beta).
\]
The degree of ${\rm deg}\,G_{\mu,\nu}$ in the space
${\cal G}_{\mu,\nu}[l_1,l_2,l_3]$ is given by
\begin{eqnarray}
D_{\mu,\nu}[l_1,l_2]&=&\sum_\alpha(\alpha-l_1)^+m_\alpha(\mu)
+\sum_\alpha(\alpha-l_2)^+m_\alpha(\nu)\nonumber\\
&+&\sum_{\alpha,\beta}A_{\alpha,\beta}m_\alpha(\mu) m_\beta(\mu)
+\sum_{\alpha,\beta}A_{\alpha,\beta}m_\alpha(\nu) m_\beta(\nu)
-\sum_{\alpha,\beta}A_{\alpha,\beta}m_\alpha(\mu) m_\beta(\nu)\label{KDEG}
\end{eqnarray}

In the special case when $l_3={\rm min}(l_1,l_2)$,
we have
\[
\tau^{(\alpha,\beta)}[l_1,l_2,{\rm min}(l_1,l_2)]\leq0.
\]
and thus there is no lower-bound condition on the riggings.

The summation (\ref{RCHAR}) 
with respect to the riggings $r,s$ can be immediately computed using
\begin{equation*}
\sum_{0\leq r_1\leq\dots\leq r_n\leq M}q^{r_1+\dots+r_n}=
\Bigl[{M+n\atop n}\Bigr],
\end{equation*}
where (if $m\in{\Z}$ and $n\in{\Z}_{\geq0}$)
the Gaussian polynomials are
\begin{eqnarray}
\Bigl[{m\atop n}\Bigr]
=\begin{cases}
\frac{\prod_{i=1}^m(1-q^i)}{\prod_{i=1}^n(1-q^i)\prod_{i=1}^{m-n}(1-q^i)}
\quad&\text{if $m\geq n$};\\
0\quad&\text{if $m<n$}.
\end{cases} 
\end{eqnarray}

\begin{lemma}
\begin{eqnarray*}
&&{\rm ch}_qW^{(M,N)}_k[l_1,l_2,{\rm min}(l_1,l_2)]_{m,n}\\
&&\quad=\sum_{|\mu|=m,|\nu|=n}q^{D_{\mu,\nu}[l_1,l_2]}
\prod_\alpha\Biggl[{P^{(M)}_{\mu,\nu}[l_1]_\al
+m_\alpha(\mu)\atop m_\alpha(\mu)}\Biggr]
\prod_\alpha\Biggl[{Q^{(N)}_{\mu,\nu}[l_2]_\al
+m_\alpha(\nu)\atop m_\alpha(\nu)}\Biggr].
\end{eqnarray*}
\end{lemma}

Let $W^{(M,N)}_k[l_1,l_2]=W^{(M,N)}_k[l_1,l_2,\min(l_1,l_2)]$. Then
\begin{theorem}\label{FERMIONIC}
The character of the space of coinvariants 
$W^{(M,N)}_k[l_1,l_2]$ is given by
\begin{equation}
\chi^{(M,N)}_k[l_1,l_2](z_1,z_2,q)
=\sum_{\mu,\nu}z_1^{|\mu|}z_2^{|\nu|}q^{D_{\mu,\nu}[l_1,l_2]}
\prod_\alpha\Biggl[{P^{(M)}_{\mu,\nu}[l_1]_\al
+m_\alpha(\mu)\atop m_\alpha(\mu)}\Biggr]
\prod_\alpha\Biggl[{Q^{(N)}_{\mu,\nu}[l_2]_\al
+m_\alpha(\nu)\atop m_\alpha(\nu)}\Biggr].
\end{equation}
\end{theorem}
In what follows, we set
 $\chi^{(M,N)}_k[l_1,l_2](z_1,z_2,q)=0$ if $l_1<0$ or $l_2<0$.

In Section 6.4 of Part II \cite{FKLMM2}, we obtained several identities
between the characters of coinvariant spaces for
$\heis$ and $\slt$-modules.
One can apply the result above to give ``fermionic'' formulas for them.
For, example, we have
\begin{theorem}
The character \Ref{sltchar} of the $\widehat{\mathfrak sl}_2$ coinvariant space
$L^{(M,N)}_{k,l}$ is
\[
\chi^{(M,N)}_{k,l}=
z^{-l}
\Bigl(\chi^{(M+1,N)}_k[l,k-l](q^{-2}z^2,z^{-2},q)
-q\chi^{(M+1,N)}_k[l-1,k-l-1](q^{-2}z^2,z^{-2},q)\Bigr).
\]
\end{theorem}

\section{The upper and lower subsets of rigged configurations}\label{RIG PREL}
In the rest of the paper we prove Theorem \ref{recursion theorem}. In this section 
we define admissible pairs $(I,J)$ of subsets of $\{1,\dots,k\}$. 
Then, we define two kinds of subsets
of rigged partitions indexed by admissible pairs, the lower and upper subsets.
and construct a bijection {}from the upper to the lower subsets indexed by
the same pair $(I,J)$.

\subsection{Admissibility of $(I,J)$}\label{ADMIS}
For a $k$-vector $\rho\in\Z^{k}$,
the $\al$-th coordinate of $\rho$ is denoted by $\rho_\al$.
The positive and negative parts of $\rho$,
$\rho^\pm\in\Z^k_{\geq 0}$, are defined by $(\rho^\pm)_\al=(\rho_\al)^\pm$.
We have $\rho=\rho^+-\rho^-$.

For $k$-vectors
$\xi,\eta\in{\Z}^k$ we write $\xi\geq\eta$ if and only if
$\xi_\alpha\geq\eta_\alpha$ for all $1\leq\alpha\leq k$. In particular,
$\xi\geq0$ means $\xi_\alpha\geq0$ for all $1\leq\alpha\leq k$.

For $I\subset\{1,\dots,k\}$, define the $k$-vectors
$\kappa(I)\in\{0,1,\dots,k\}^k$,  $\varepsilon(I)\in\{-1,0,1\}^k$ by the formula
\begin{align*}
\kappa(I)_\al&=\sum_{i\in I,\;\;i\leq \al} 1,\\
\varepsilon(I)_\al&=\sum_{i\in I}(\delta_{i,\al}-\delta_{i,\al+1}),
\end{align*}
where $\al=1,\dots,k$. 

We define a partial ordering in the set $2^{\{1,\dots,k\}}$:
$J\geq J'$ if and only if $\kappa(J)\geq\kappa(J')$.
If we set $J=\{v_1,\dots,v_s\}$ and $J'=\{v'_1,\dots,v'_{s'}\}$
where $v_1<\dots<v_s$ and $v'_1<\dots<v'_{s'}$, this is equivalent to
$s\geq s'$ and $v_i\leq v'_i$ for $1\leq i\leq s$.

Sometimes it is convenient to extend the definition of $\kappa(I)$
to $I$ not necessarily satisfying $I\subset\{1,\dots,k\}$.
Namely, we use the same definition for $I\subset\{1,2,3,\dots\}$.
Note, however, that $\kappa(I)_\alpha=\kappa(I\cap\{1,\dots,k\})_\alpha$
because we consider $\alpha$ only in the region $\{1,\dots,k\}$.

Note that if $I=I_1\bigsqcup I_2$ then $\kappa(I)=\kappa(I_1)+\kappa(I_2)$
and $\varepsilon(I)=\varepsilon(I_1)+\varepsilon(I_2)$. We have
$\kappa(I)=\sum_{i\in I}\kappa(i)$ and $\varepsilon(I)=\sum_{i\in I}\varepsilon(i)$,
where we denoted  $\kappa(i)=\kappa(\{i\})$ and $\varepsilon(i)=\varepsilon(\{i\})$.
For example, if $k=5$, we have  $\kappa(\{2,4,5\})=(0,1,1,2,3)$ and
$\varepsilon(\{2,4,5\})=(-1,1-1,0,1)$. 
For $\al,\beta\in\{1,2,3,\dots\}$, we denote the interval
$\{\alpha,\alpha+1,\dots,\beta\}$ by $[\alpha,\beta]$ and
the $k$-vector $\kappa([\alpha,\beta])$ by $\kappa[\alpha,\beta]$.

Fix $0\leq l_1,l_2\leq k$.
Let $I=\{u_1,\dots,u_a\}$ ($u_1<\dots<u_a$) and $J=\{v_1,\dots,v_b\}$
($v_1<\dots <v_b$) be subsets of $\{1,\dots,k\}$.
We define the $(l_1,l_2)$-admissibility of $(I,J)$ as follows.

Let $p=p(l_1,J)$ be the number of elements of $J$ which are less than $l_1+1$,
\be
v_1<\dots<v_p<l_1+1\leq v_{p+1}<\dots <v_b.
\ee

Set
\begin{equation}
t={\rm max}(1,l_1+b-k+1).
\end{equation}
We have $t\leq a$ if and only if $l_1+c<k$.
\begin{lemma}
We can label the complement of $[l_1+1,k]\cap J$ in $[l_1+1,k]$ as follows.
\begin{eqnarray}
&&[l_1+1,k]\backslash\{v_{p+1},\dots,v_b\}\nonumber\\
&&=
\begin{cases}
\{v'_p,\dots,v'_t\}\quad\text{where $v'_p<\dots<v'_t$}
&\text{if $l_1+b\geq k$};\\
\{v'_p,\dots,v'_1,w_1,\dots,w_{k-l_1-b}\}
\quad\text{where $v'_p<\dots<v'_1<w_1<\dots<w_{k-l_1-b}$}
&\text{if $l_1+b<k$}.
\end{cases}\nonumber\\
\end{eqnarray}
\end{lemma}
\begin{proof}
Note that $\card([l_1+1,k]\backslash\{v_{p+1},\dots,v_b\})=k-l_1-b+p$.
If $l_1+b\geq k$ we have $k-l_1-b+p\leq p$, and we label the complement as
$v'_p<\dots<v'_t$. If $l_1+b<k$ we have $k-l_1-b+p>p$, and
we label the complement as $v'_p<\dots<v'_1<w_1<\dots<w_{k-l_1-b}$.
\end{proof}

In the case $l_1+b\geq k$ it is convenient to set 
\begin{equation}
v'_{t-1}=\dots=v'_1=k+1.\label{CONV}
\end{equation}

We have
\begin{lemma}
\bea\label{prime definition}
\Bigl(\kappa(J)-\kappa[l_1+1,l_1+b)]\Bigr)^+
=\sum_{i=1}^p\left(\kappa(v_i)-\kappa(v_i')\right).
\ena 
\end{lemma}
\begin{proof}
Observe $\kappa_\alpha=(\kappa(J)-\kappa[l_1+1,l_1+b])_\alpha$
when $\alpha$ varies {}from $1$ to $k$.
For $\alpha\leq l_1$, $\kappa$ increases {}from $\kappa_{\alpha-1}$
to $\kappa_\alpha$ by $1$ if $\alpha=v_i$
($1\leq i\leq p$) or stays constant otherwise. In the case $l_1+b\geq k$,
for $l_1+1\leq\alpha\leq k$, $\kappa$ decreases {}from $\kappa_{\alpha-1}$
to $\kappa_\alpha$ by $1$ if $\alpha=v'_i$
($t\leq i\leq p$) or stays constant otherwise. In the case $l_1+b<k$,
for $l_1+1\leq\alpha\leq v'_1$, $\kappa_\alpha$ decreases by $1$
at $\alpha=v'_i$ ($1\leq i\leq p$) or stays constant otherwise. 
In particular, we have $\kappa_{v'_1}=0$. For $\alpha>v'_1$,
$\kappa_\alpha\leq0$. The equality (\ref{prime definition})
follows {}from these observations
with the convention (\ref{CONV}) for $l_1+b\geq k$.
\end{proof}

A pair of subsets $(I,J)$ is called $(l_1,l_2)$-admissible if
\bea\label{admissible}
a\leq p(l_1,J),\quad b\leq l_2\text\quad{\rm and}\quad v_i\leq u_i<v_i'\,\,(1\leq i\leq a).
\ena
Note that $l_2$ appears only in the restriction $b\leq l_2$.
A pair $(\emptyset,J)$ is $(l_1,l_2)$-admissible if and only if
$\card(J)\leq l_2$. An $(l_1,k)$-admissible pair is simply called
$l_1$-admissible. If $(I,J)$ is $l_1$-admissible, then $\card(I)\leq l_1$.

Note that if $l_1+c\geq k$, $(I,J)$ is $l_1$ admissible if and only if
\[
v_i\leq u_i\quad(1\leq i\leq a).
\]
The condition $a\leq p(l_1,J)$ is satisfied because
$v_a\leq v_b-c\leq k-c\leq l_1$.

If $l_1+c<k$, for an $l_1$-admissible pair $(I,J)$ we set
\begin{eqnarray}
&&\tilde I=\tilde I(I,J)=I\sqcup I'\label{DEFTI}\\
&&I'=
\begin{cases}
\{v'_a,\dots,v'_t\}\quad&\text{if $l_1+b\geq k$};\\
\{v'_a,\dots,v'_1\}\sqcup\{w_1,\dots,w_{k-l_1-b}\}
\quad&\text{if $l_1+b<k$}.
\end{cases}
\end{eqnarray}
Note that
\begin{equation}\label{IEQ}
\card(\tilde I)=k-l'_1
\end{equation}
We also set
\begin{equation}
\tilde J=J\cap[1,v'_a-1].\label{DEFTJ}
\end{equation}
We have
\begin{equation}\label{JEQ}
\card(\tilde J)=v'_a+c-l'_1-1
\end{equation}
because
\begin{eqnarray*}
\card(\tilde J)&=&\card(J\cap[1,l_1])+\card(J\cap[l_1+1,v'_a-1])\\
&=&p+\card([l_1+1,v'_a-1])-\card(\{v'_p,\dots,v'_{a+1}\})\\
&=&a+v'_a-l_1-1\\
&=&v'_a+c-l'_1-1.
\end{eqnarray*}

If $l_1+c\geq k$ we have $v'_a=k+1$ by (\ref{CONV}). We set
$\tilde I=I$ and $\tilde J=J$. The equalities (\ref{IEQ}) and (\ref{JEQ})
are valid in this case, too.

\begin{lemma}
Suppose that $l_1+c<k$. The map
${\mathfrak b}:(I,J)\mapsto(\tilde I,\tilde J)$ given by (\ref{DEFTI})
and (\ref{DEFTJ}) is a bijection between the set of $l_1$-admissible pairs
$(I,J)$ satisfying $\card(I)=a$ and $\card(J)=a+c$ and the set of
$(\tilde I,\tilde J)$ satisfying
\begin{eqnarray}
&&\tilde I=\{u_1,\dots,u_{\tilde a}\}\,\,(\tilde a=a+k-l_1-c),
\quad u_1<\dots<u_{\tilde a},
\quad u_{a+1}\geq l_1+1,\label{TIL-I}\\
&&\tilde J=\{v_1,\dots,v_{\tilde b}\}\,\,(\tilde b=a+u_{a+1}-l_1-1),
\quad v_1<\dots<v_{\tilde b},\notag\\
&&v_a\leq l_1,\quad v_{\tilde b}<u_{a+1},
\quad v_i\leq u_i\quad(1\leq i\leq a).
\label{TIL-J}
\end{eqnarray}
$($In $(\ref{TIL-J})$, the condition $v_a\leq l_1$ follows {}from the others.$)$
\end{lemma}
\begin{proof}
We will prove that the inverse map
${\mathfrak c}:(\tilde I,\tilde J)\mapsto(I,J)$ is given by
\begin{equation}
I=\{u_1,\dots,u_a\},\quad
J=\tilde J\sqcup([u_{a+1},k]\backslash\tilde I).
\label{IJ}
\end{equation}

Let us prove that the composition ${\mathfrak c}\circ{\mathfrak b}$
is the identity map. Consider $(I,J)$ and
$(\tilde I,\tilde J)={\mathfrak b}(I,J)$.
Set $(I_1,J_1)={\mathfrak c}(\tilde I,\tilde J)$.
Since $u_a<v'_a$, the smallest $a$ elements in $\tilde I$ are $u_1<\dots<u_a$.
Therefore, $I=I_1$.

Note that 
\begin{eqnarray*}
v'_a&=&u_{a+1},\\
\tilde J&=&J\cap[1,u_{a+1}-1],\\
\text{$[u_{a+1},k]$}\backslash\tilde I
&=&\text{$[u_{a+1},k]$}\backslash I'\\
&=&[u_{a+1},k]\backslash\Bigl([u_{a+1},k]\cap([l_1+1,k]\backslash J)\Bigr)\\
&=&[u_{a+1},k]\backslash([u_{a+1},k]\backslash J)\\
&=&[u_{a+1},k]\cap J.
\end{eqnarray*}
Therefore, we have
\[
J_1=\tilde J\sqcup([u_{a+1},k]\backslash\tilde I)=
(J\cap[1,u_{a+1}-1])\sqcup([u_{a+1},k]\cap J)
=J.
\]

Let us prove that the pair $(I,J)$ given by (\ref{IJ})
is $l_1$-admissible. We define $u_i$ and $v_j$ as before {}from $I$ and $J$.
It is clear that $v_i\leq u_i$ ($1\leq i\leq a$).

The number $p=p(l_1,J)$ satisfies $v_{p+1}>l_1$. Since $v_a\leq l_1$,
we have $a\leq p$. 

Let the smallest $p-a+1$ elements of the set $[l_1+1,k]\backslash J$ be
$\{v'_p,\dots,v'_a\}$ ($v'_p<\dots<v'_a$). We will show that  $v'_a=u_{a+1}$.
Then, it follows that $u_i<v'_i$ ($1\leq i\leq a$).

Since $v_{\tilde b}<u_{a+1}$, we have $\tilde J\cap[u_{a+1},k]=\emptyset$.
Then, we have
\begin{eqnarray}
[l_1+1,k]\backslash J&=&
([l_1+1,u_{a+1}-1]\backslash\tilde J)\sqcup
\Bigl([u_{a+1},k]\backslash([u_{a+1},k]\backslash\tilde I)\Bigr)\nonumber\\
&=&([l_1+1,u_{a+1}-1]\backslash\tilde J)\sqcup
([u_{a+1},k]\cap\tilde I).
\end{eqnarray}
Since
$\card([l_1+1,u_{a+1}-1]\backslash\tilde J)=u_{a+1}-1-l_1-(\tilde b-p)=p-a$,
we have $v'_a=u_{a+1}$.
\end{proof}

Fix $a,c$ and $I$ with $\card(I)=a$.
In Section \ref{upper summation section} we will use
the minimal element $J_{\rm min}$ among $J$ such that
$\card(J)=a+c$, $(I,J)$ is $l_1$-admissible and $\tilde I(I,J)$ is fixed.
\begin{lemma}\label{MIN-1}
Suppose that $l_1+c\geq k$ and fix $I=\{u_1,\dots,u_a\}$.
Consider the set of $J$ such that
$(I,J)$ is $l_1$-admissible and $\card(J)=a+c$.
This set has the minimal element given by
\begin{equation}\label{MIN1}
J_{\rm min}=\{{\rm min}(u_i,k-b+i)\}_{1\leq i\leq a}\sqcup[k-c+1,k].
\end{equation}
\end{lemma}
\begin{proof}
Set $J=\{v_1,\dots,v_b\}$ ($v_1<\dots<v_b$). We have obviously
$v_i\leq k-b+i$ ($1\leq i\leq b$). For $1\leq i\leq a$ we have further
$v_i\leq u_i$. Therefore, the minimal element is given by (\ref{MIN1}).
\end{proof}
For $l_1+c<k$, we obtain $J_{\rm min}$ by using the bijection ${\mathfrak b}$.
\begin{lemma}\label{MIN-2}
Suppose that $l_1+c<k$ and fix $\tilde I$ satisfying $(\ref{TIL-I})$.
Set $I=\{u_1,\dots,u_a\}$ and consider the set of $J$ such that
$(I,J)$ is $l_1$-admissible, $\card(J)=a+c$ and $\tilde I(I,J)=\tilde I$.
This set has the minimal element given by
\begin{equation}
J_{\rm min}=\{{\rm min}(u_i,l_1-a+i)\}_{1\leq i\leq a}\sqcup[l_1+1,u_{a+1}-1]
\sqcup([u_{a+1},k]\backslash\tilde I).\label{MIN2}
\end{equation}
\end{lemma}
\begin{proof}
The minimal set $\tilde J_{\rm min}$ among $\tilde J$ satisfying (\ref{TIL-J})
is given by
$\tilde J_{\rm min}
=\{{\rm min}(u_i,l_1-a+i)\}_{1\leq i\leq a}\sqcup[l_1+1,u_{a+1}-1]$.
Then, $J_{\rm min}$ is given by (\ref{IJ}).
\end{proof}
\subsection{Vectors $\rho$ and $\sigma$ and lower subsets}\label{LOWSEC}
For an $(l_1,l_2)$-admissible pair $(I,J)$, define the vectors
$\rho(I,J)=\rho_{l_1,k}(I,J)$, $\sigma(J)=\sigma_{l_2,k}(J)\in\Z^k$ by
\begin{eqnarray}
\rho(I,J)&=&\sum_{i=1}^a(\kappa(v_i)-\kappa(u_i))
+\sum_{i=a+1}^p(\kappa(v_i)-\kappa(v_i')),\label{REALRHO}\\
\sigma(J)&=&\kappa[1,l_2]-\kappa(J).\label{REALSIGMA}
\end{eqnarray}
Note that $\rho(I,J),\sigma(J)\geq0$.

We introduce a few notations.

We use the symbols
\begin{eqnarray}
\leq_\varepsilon&\leftrightarrow&
\begin{cases}
\leq&\text{ if $\varepsilon\not=1$;}\label{LOWARROW}\\
=&\text{ if $\varepsilon=1$,}\\
\end{cases}\\
\leq^\varepsilon&\leftrightarrow&
\begin{cases}
\leq&\text{ if $\varepsilon\not=-1$;}\\
=&\text{ if $\varepsilon=-1$.}
\end{cases}\label{UPARROW}
\end{eqnarray}

For a rigging
$r=\{r_i^{(\al)}\}_{1\leq\al\leq k\atop1\leq i\leq m_\al}$, we define
\begin{eqnarray}
r[\alpha]=
\begin{cases}
r^{(\alpha)}_{m_\alpha}&\text{ if $m_\alpha\geq1$;}\\
\infty&\text{ if $m_\alpha=0$.}\label{INF}
\end{cases}
\end{eqnarray}

For a pair of subsets $(I,J)$ ($I,J\subset\{1,\dots,k\}$)
and a pair of integers $(l_1,l_2)$ ($0\leq l_1,l_2\leq k$), we define
the subset $R_{m,n}[l_1,l_2]_{I,J}\subset R_{m,n}$ as follows. If
$(I,J)$ is $(l_1,l_2)$-admissible, we set
\begin{eqnarray}
&&R_{m,n}[l_1,l_2]_{I,J}=
\{(\mu,r;\nu,s)\in R_{m,n};\nonumber\\
&&\quad
r[\alpha]\geq_{\varepsilon(I)_\alpha}\rho(I,J)_\al\quad\text{\rm and}\quad
s[\alpha]\geq_{\varepsilon(J)_\alpha}\sigma(J)_\al
\quad\text{\rm for all $1\leq\alpha\leq k$}\}
\end{eqnarray}
(see (\ref{LOWARROW})). Otherwise we set 
$R_{m,n}[l_1,l_2]_{I,J}=\emptyset$.

The restriction for $r[\alpha]$ is called marked if
$\varepsilon(I)_\alpha=1$ and, therefore, it takes the form
$\rho(I,J)_\alpha=r[\alpha]$; it is called unmarked otherwise, namely,
if it takes the form $\rho(I,J)_\alpha\leq r[\alpha]$. Similarly,
we distinguish the marked and unmarked restrictions for $s[\alpha]$.

Suppose that $(\mu,r;\nu,s)$ is contained in $R_{m,n}[l_1,l_2]_{I,J}$.
Then, $m_\alpha\not=0$ if $\varepsilon(I)_\alpha=1$;
$n_\alpha\not=0$ if $\varepsilon(J)_\alpha=1$.

For $M,N\geq0$ we define
\begin{equation}
R^{(M,N)}_{m,n}[l_1,l_2]_{I,J}=R_{m,n}[l_1,l_2]_{I,J}\cap
R^{(M,N)}_{m,n}[l_1,l_2].
\end{equation}

If an element $(\mu,r;\nu,s)$ is contained in
$R^{(M,N)}_{m,n}[l_1,l_2]_{I,J}$, and if $m_\alpha\not=0$ for some $\alpha$,
then we have $\rho(I,J)_\alpha\leq P^{(M)}_{\mu,\nu}[l_1]_\alpha$; if
$n_\alpha\not=0$ then $\sigma(J)_\alpha\leq Q^{(N)}_{\mu,\nu}[l_2]_\alpha$.
In the rest of this section we prove the validity of these inequalities
when $m_\alpha=0$ or $n_\alpha=0$.

Let us abbreviate $P^{(M)}_{\mu,\nu}[l_1]_\alpha$ to $P_\alpha$,
and $\rho(I,J)_\alpha$ to $\rho_\alpha$.
Recall that
\begin{equation}
\rho(I,J)=\sum_{\alpha=1}^p\kappa(v_\alpha)
-\sum_{\alpha=1}^a\kappa(u_\alpha)-\sum_{\alpha=a+1}^p\kappa(v'_\alpha),
\label{CUT0}
\end{equation}
where $v_1<\dots<v_p\leq l_1$, $v'_p<\dots<v'_t<v'_{t-1}=\dots=v'_1=k+1$,
$u_1<\dots<u_a$
and $v_i\leq u_i< v'_i$ ($1\leq i\leq a$). Here $t={\rm max}(1,l_1+b-k+1)$.

We set
\begin{eqnarray}
J_{\rm up}&=&\{v_1,\dots,v_p\}=J\cap[1,l_1],\label{JUP}\\
J'_{\rm down}&=&
\begin{cases}
\{v'_p,\dots,v'_t\}&\text{ if $l_1+c\geq k$;}\\
\{v'_p,\dots,v'_{a+1}\}&\text{ if $l_1+c<k$.}
\end{cases}
\end{eqnarray}

If $l_1+c\geq k$ we have $t>a$ and
\begin{equation}
J'_{\rm down}=[l_1+1,k]\backslash J.\label{DOWN}
\end{equation}

If $l_1+c<k$ we have $t\leq a$ and
\begin{equation}
J'_{\rm down}\cup
\Bigl(\{\alpha;\alpha\geq l_1+1,\rho_\alpha=0\}\backslash J\Bigr)
=[l_1+1,k]\backslash J.\label{DOWN2}
\end{equation}

We list a few more properties of $\rho_\alpha$.
\begin{eqnarray}
&(P1)&\quad\rho_\alpha-2\leq\rho_{\alpha+1}\leq\rho_\alpha+1,\nonumber\\
&(P2)&\quad\text{If $\alpha+1\in J_{\rm up}$,
then $\rho_{\alpha+1}\geq\rho_\alpha$,}\nonumber\\
&(P3)&\quad\text{If $\alpha+1\not\in J_{\rm up}$,
then $\rho_{\alpha+1}\leq\rho_\alpha$,}\nonumber\\
&(P4)&\quad\text{If $\alpha+1\in J'_{\rm down}$,
then $\rho_{\alpha+1}\leq\rho_\alpha-1$,}\nonumber\\
&(P5)&\quad\text{If $\alpha+1\not\in J'_{\rm down}$,
then $\rho_{\alpha+1}\geq\rho_\alpha-1$.}\nonumber
\end{eqnarray}
\begin{lemma}
If $R^{(M,N)}_{m,n}[l_1,l_2]_{I,J}$ contains an element $(\mu,r;\nu,s)$,
then we have $\rho_k\leq P_k.$.
\end{lemma}
\begin{proof}
Assume that $P_i\geq\rho_i$ and $P_\alpha<\rho_\alpha$ ($i+1\leq\alpha\leq k$)
for some $0\leq i\leq k-1$. As we noted at the beginning of this section
we have $m_\alpha=0$ ($i+1\leq \alpha\leq k$). This implies
$\varepsilon(I)_\alpha\not=1$ ($i+1\leq \alpha\leq k$). Therefore, we have
\begin{equation}
I\subset[1,i].\label{IRES}
\end{equation}

We have
\begin{eqnarray}
P_{i+1}-P_i&=&M-(i+1-l_1)^++(i-l_1)^+
+\sum_{\beta\geq i+1}n_\beta\nonumber\\
&<&\rho_{i+1}-\rho_i.\label{CUT3}
\end{eqnarray}

\noindent
{\it Subcase 1} : $i+1\leq l_1$.

{}{}From $(P1)$ we have $\rho_{i+1}-\rho_i\leq1$. Using (\ref{CUT3}) we have
$$
P_{i+1}-P_i=M+\sum_{\beta\geq i+1}n_\beta<1.
$$
Therefore  we have $M=0$ and $n_{i+1}=\dots=n_k=0$. This implies
$J\subset[1,i]$, and therefore $i+1\not\in J$.
Using $(P3)$ we have $\rho_{i+1}\leq\rho_i$. This is a contradiction because
$$
0=P_{i+1}-P_i<\rho_{i+1}-\rho_i\leq0.
$$

\noindent
{\it Subcase 2} : $i+1\geq l_1+1$.

We have $i+1\not\in J_{\rm up}$ because $J_{\rm up}\subset[1,l_1]$.
{}{}From $(P3)$ follows $\rho_{i+1}-\rho_i\leq0$ and using (\ref{CUT3}) we have
$$
P_{i+1}-P_i=M-1+\sum_{\beta\geq i+1}n_\beta<0.
$$
Therefore, we have $M=0$, $n_{i+1}=\dots=n_k=0$ and $i+1\not\in J$ again.

If $l_1+c\geq k$, because of (\ref{DOWN}) we have $i+1\in J'_{\rm down}$.
Using $(P4)$ we have $\rho_{i+1}-\rho_i\leq-1$. This is a contradiction because
$$
-1=P_{i+1}-P_i<\rho_{i+1}-\rho_i\leq-1.
$$

If $l_1+c<k$, we proceed as follows. If $i+1\in J'_{\rm down}$,
it leads to a contradiction as above. If $i+1\not\in J'_{\rm down}$,
because of (\ref{DOWN2}) we have $\rho_{i+1}=0$.
It implies $P_{i+1}<0$. However, this is prohibited by (\ref{COC}).
\end{proof}
\begin{lemma}\label{REFCOM}
If $R^{(M,N)}_{m,n}[l_1,l_2]_{I,J}$ contains an element $(\mu,r;\nu,s)$,
then we have $\rho_\alpha\leq P_\alpha$ ($1\leq \alpha<k$).
\end{lemma}
\begin{proof}
We lead to a contradiction assuming that
for some $i$ and $j$ satisfying $1\leq i+1<j\leq k$ we have
$$
P_i\geq\rho_i,\quad P_\alpha<\rho_\alpha\,\,(i+1\leq\alpha\leq j-1),
\quad P_j\geq\rho_j.
$$

We set $p=\frac1{j-i}$ so that $pi+(1-p)j=j-1$. A simple calculation
as (\ref{REFCOM2}) shows
\begin{eqnarray}
P_{j-1}-\rho_{j-1}&\geq&\Delta\rho
+\theta(i<l_1<j)(l_1-i)p+\sum_{i<\beta<j}(\beta-i)pn_\beta,
\label{ESTIMATE}\\
\Delta\rho&=&p\rho_i+(1-p)\rho_j-\rho_{j-1}.
\end{eqnarray}
Here we used $m_\beta=0$ for $i<\beta<j$.
Note that the last two terms in the RHS of (\ref{ESTIMATE}) is non-negative.

We consider three cases
$\rho_j\geq\rho_{j-1}$, $\rho_j=\rho_{j-1}-1$ and $\rho_j=\rho_{j-1}-2$,
separately.

\noindent
{\it Case 1} : $\rho_j\geq\rho_{j-1}$.

Because of $(P1)$ we have
$\rho_{j-1}\leq\rho_i+j-i-1$. {}From this follows
$$
\Delta\rho\geq p(\rho_{j-1}-(j-i-1))+(1-p)\rho_{j-1}-\rho_{j-1}=-1+p.
$$
Using (\ref{ESTIMATE}) we have $P_{j-1}-\rho_{j-1}\geq0$, which is a
contradiction.

\noindent
{\it Case 2} : $\rho_j=\rho_{j-1}-1$.

\noindent
{\it Subcase 1} : $i\geq l_1$.

Using $(P3)$ and (\ref{JUP}) we have $\rho_{j-1}\leq\rho_i$. Then, we have
\begin{equation}
P_{j-1}-\rho_{j-1}\geq\Delta\rho\geq p\rho_{j-1}+(1-p)(\rho_{j-1}-1)-\rho_{j-1}
=-1+p.
\label{CONTRA}
\end{equation}
This is a contradiction.

\noindent
{\it Subcase 2} : $i<l_1<j$.

Because of (P3) and (\ref{JUP}) we have
$\rho_{j-1}\leq \rho_i+l_1-i$. Therefore, noting that $\theta(i<l_1<j)=1$,
we have again
\begin{eqnarray}
P_{j-1}-\rho_{j-1}
&\geq&p(\rho_{j-1}-(l_1-i))+(1-p)(\rho_{j-1}-1)-\rho_{j-1}+(l_1-i)p
\nonumber\\
&=&-1+p.\nonumber
\end{eqnarray}
This is a contradiction.

\noindent
{\it Subcase 3} : $j\leq l_1$.

We have $j\not\in J$ because otherwise
$j\in J_{\rm up}$ and and using $(P2)$ we have $\rho_j\geq\rho_{j-1}$,
which is a contradiction.

We will prove by induction the following statements for $i+1\leq\alpha\leq j-1$:
\begin{eqnarray}
&(C1)_\alpha&\quad n_\alpha=0,\nonumber\\
&(C2)_\alpha&\quad[\alpha,j]\cap J=\emptyset,\nonumber\\
&(C3)_\alpha&\quad\rho_i\geq\rho_{j-1}-(\alpha-i-1).\nonumber
\end{eqnarray}
Then, $(C3)_{i+1}$ leads to (\ref{CONTRA}), which is a contradiction.

We first note that $(C2)_j$ and $(C3)_j$ are valid. These are the basis
for the induction. {}From $(C1)_\alpha$ follows $\varepsilon(J)_\alpha\not=1$.
Using $(C2)_{\alpha+1}$ we have $\alpha\not\in J$, and therefore $(C2)_\alpha$.
Because of $(P3)$ {}from $(C2)_\alpha$ follows $(C3)_\alpha$. Finally,
we show that for $i+2\leq\alpha\leq j$ {}from $(C3)_\alpha$ follows
$(C1)_{\alpha-1}$. Unless $n_{\alpha-1}=0$ we have again
\begin{eqnarray}
P_{j-1}-\rho_{j-1}&\geq&
p(\rho_{j-1}-(\alpha-i-1))+(1-p)(\rho_{j-1}-1)-\rho_{j-1}
+(\alpha-1-i)pn_{\alpha-1}\nonumber\\
&\geq&-1+p.\nonumber
\end{eqnarray}

\noindent
{\it Case 3} : $\rho_j=\rho_{j-1}-2$.

\noindent
{\it Subcase 1} : $i\geq l_1$.

We will prove by induction the following statements for
$i+1\leq\alpha\leq j-1$.
\begin{eqnarray*}
&(C1)'_\alpha&\quad n_\alpha=0,\\
&(C2)'_\alpha&\quad[\alpha,j]\subset J'_{\rm down}\\
&(C3)'_\alpha&\quad\rho_i\geq\rho_{j-1}+j-\alpha.
\end{eqnarray*}
Then, {}from $(C3)'_{i+1}$, we have
$\rho_i\geq\rho_{j-1}+j-i-1$. Using this we have
\begin{eqnarray*}
P_{j-1}-\rho_{j-1}&\geq&p(\rho_{j-1}+j-i-1)+(1-p)(\rho_{j-1}-2)
-\rho_{j-1}\\
&=&-1+p.
\end{eqnarray*}
This is a contradiction.

As we have noted above we have $(C2)'_j$.
Because of $(P3)$ and (\ref{JUP}),
we have $\rho_i\geq\rho_{j-1}$. This is $(C3)'_j$.

Assume that $(C1)'_\alpha$ and $(C2)'_{\alpha+1}$ are valid
for some $i+1\leq\alpha\leq j-1$.
{}{}From $(C1)'_\alpha$ follows $\varepsilon(J)_\alpha\not=1$. Since
$\alpha+1\not\in J$ by $(C2)'_{\alpha+1}$, we have $\alpha\not\in J$.

If $l_1+c\geq k$, because of (\ref{DOWN}) we have $\alpha\in J'_{\rm down}$.
If $l_1+c<k$, we use (\ref{DOWN2}). Note that $\alpha\geq i+1\geq l_1+1$
and $\alpha\not\in J$. If $\rho_\alpha=0$,
we have $P_\alpha<0$, which contradicts (\ref{COC}).
Otherwise, we have $\alpha\in J'_{\rm down}$.
Thus we have derived $(C2)'_\alpha$ {}from $(C1)'_\alpha$ and $(C2)'_{\alpha+1}$.

Using (P3) and (P4) we can derive $(C3)'_\alpha$ {}from $(C2)'_\alpha$ .

Suppose that we have $(C3)'_\alpha$ for some $i+2\leq\alpha\leq j$.
Unless $n_{\alpha-1}=0$ we have
\begin{eqnarray*}
P_{j-1}-\rho_{j-1}&\geq&p(\rho_{j-1}+j-\alpha)+(1-p)(\rho_{j-1}-2)
-\rho_{j-1}+(\alpha-i-1)pn_{\alpha-1}\\
&\geq&-1+p.
\end{eqnarray*}
This is a contradiction. We have derived $(C1)'_{\alpha-1}$
{}from $(C3)'_\alpha$.

\noindent
{\it Subcase 2} : $i<l_1<j$. 

We will prove by induction the following statements for
$i+1\leq\alpha\leq j-1$.
\begin{eqnarray*}
&(C1)''_\alpha&\quad n_\alpha=0,\\
&(C2)''_\alpha&\quad[\alpha,j]\cap J=\emptyset,\\
&(C3)''_\alpha&\quad\rho_i\geq\rho_{j-1}-l_1-\alpha+i+j.
\end{eqnarray*}
Then, {}from $(C3)''_{i+1}$ we have
$\rho_i\geq\rho_{j-1}+j-l_1-1$. Therefore we have
$$
P_{j-1}-\rho_{j-1}\geq p(\rho_{j-1}+j-l_1-1)
+(1-p)(\rho_{j-1}-2)-\rho_{j-1}+(l_1-i)p=-1+p,
$$
which is a contradiction.

We have $(C2)''_j$ and $(C3)''_j$. It is obvious that {}from $(C1)''_\alpha$
and $(C2)''_{\alpha+1}$ follows $(C2)''_\alpha$.

Suppose that $(C2)''_\alpha$ is valid for some $i+1\leq\alpha\leq j-1$.
In particular, we have $\alpha\not\in J$.
If $\alpha\geq l_1+1$, using (\ref{DOWN}) or (\ref{DOWN2}) we have
$\alpha\in J'_{\rm down}$ unless we have $l_1+c<k$ and
$P_\alpha<\rho_\alpha=0$, which contradicts (\ref{COC}).
Therefore, by using (P4) (if $\alpha\geq l_1+1$)
or (P3) (if $\alpha\leq l_1$) we have $(C3)''_\alpha$.

Suppose that $(C3)''_\alpha$ is valid for some $i+2\leq\alpha\leq j$.
Unless $n_{\alpha-1}=0$ {}from (\ref{ESTIMATE}) we have
\begin{eqnarray*}
P_{j-1}-\rho_{j-1}&\geq&p(\rho_{j-1}-l_1-\alpha+i+j))
+(1-p)(\rho_{j-1}-2)-\rho_{j-1}\\
&&\,\,+(l_1-i)p+(\alpha-i-1)pn_{\alpha-1}\\
&\geq&-1+p.
\end{eqnarray*}
This is a contradiction. Thus, we have proved $(C1)''_{\alpha-1}$.

\noindent
{\it Subcase 3} : $l_1\geq j$.

Because of $(P5)$ we have $j\in J'_{\rm down}$.
Because of (\ref{DOWN}), this is a contradiction.
\end{proof}

Next we proceed to the inequality $\sigma(J)\leq Q_{\mu,\nu}^{(N)}[l_2]$.
Let us abbreviate 
$Q_{\mu,\nu}^{(N)}[l_2]_\alpha$ to $Q_\alpha$ and
$\sigma(J)_\alpha$ to $\sigma_\alpha$. Recall that $b\leq l_2$ and
\begin{equation}
\sigma(J)=\sum_{\alpha=1}^{l_2}\kappa(\alpha)-\sum_{i=1}^b\kappa(v_i).
\label{RB}
\end{equation}
We have, in particular,
$\sigma_\alpha-1\leq\sigma_{\alpha+1}\leq\sigma_\alpha+1$.
\begin{lemma}
Suppose that $N\geq1$.
If $R^{(M,N)}_{m,n}[l_1,l_2]_{I,J}$ contains an element $(\mu,r;\nu,s)$,
then we have $\sigma_k\leq Q_k$.
\end{lemma}
\begin{proof}
Assume that $Q_i\geq\sigma_i$ for some $0\leq i\leq k-1$
and $Q_\alpha<\sigma_\alpha$ ($i+1\leq\alpha\leq k$).
We have $n_\alpha=0$ ($i+1\leq\alpha\leq k$). This implies
$\varepsilon(J)_\alpha\not=1$ ($i+1\leq\alpha\leq k$). Therefore, we have
$J\subset[1,i]$.

We have
\begin{eqnarray}
Q_{i+1}-Q_i&=&N-(i+1-l_2)^++(i-l_2)^++\sum_{\beta\geq i+1}m_\beta\notag\\
&<&\sigma_{i+1}-\sigma_i.\label{RC}
\end{eqnarray}

\noindent
{\it Subcase 1} : $i+1\leq l_2$.

{}{}From (\ref{RB}) we have $\sigma_{i+1}-\sigma_i\leq1$. Using (\ref{RC}) we have
$$
Q_{i+1}-Q_i=N+\sum_{\beta\geq i+1}m_\beta<1.
$$
This is a contradiction because we assumed $N-1\geq0$.

\noindent
{\it Subcase 2} : $l_2\leq i$.

We have
$$
Q_{i+1}-Q_i=N-1+\sum_{\beta\geq i+1}m_\beta<0.
$$
This is a contradiction.
\end{proof}
\begin{lemma}
Suppose that $N\geq1$.
If $R^{(M,N)}_{m,n}[l_1,l_2]_{I,J}$ contains an element $(\mu,r;\nu,s)$,
then we have $\sigma_\alpha\leq Q_\alpha$ ($1\leq \alpha<k$).
\end{lemma}
\begin{proof}
Suppose that for some $i$ and $j$ such that $1\leq i+1<j\leq k$ we have
$Q_i\geq\sigma_i$, $Q_\alpha<\sigma_\alpha$ ($i+1\leq\alpha\leq j-1$)
and $Q_j\geq\sigma_j$. We set $p=\frac1{j-i}$. We have
\begin{eqnarray}
Q_{j-1}-\sigma_{j-1}&\geq&\Delta\sigma+\theta(i<l_2<j)(l_2-i)p
+\sum_{i<\beta<j}(\beta-i)pm_\beta,\label{RD}\\
\Delta\sigma&=&p\sigma_i+(1-p)\sigma_j-\sigma_{j-1}.
\end{eqnarray}

\noindent
{\it Case 1} : $\sigma_j\geq\sigma_{j-1}$.

We have
$$
\Delta\sigma\geq p\Bigl(\sigma_{j-1}-(j-i-1)\Bigr)
+(1-p)\sigma_{j-1}-\sigma_{j-1}=-1+p.
$$
Using (\ref{RD}) we have $Q_{j-1}-\sigma_{j-1}\geq0$, which is a contradiction.

\noindent
{\it Case 2} : $\sigma_j=\sigma_{j-1}-1$. 

{}{}From (\ref{RB}) we have $l_2+1\leq j$ and
$\sigma_{j-1}\leq\sigma_j+(l_2-i)^+$.
Therefore, we have
$$
Q_{j-1}-\sigma_{j-1}\geq p\Bigl(\sigma_{j-1}-(l_2-i)^+\Bigr)
+(1-p)(\sigma_{j-1}-1)-\sigma_{j-1}+(l_2-i)^+p=-1+p,
$$
which is a contradiction.
\end{proof}

We have proved
\begin{proposition}\label{NON}
Suppose that $M,N-1\geq0$. If $R^{(M,N)}_{m,n}[l_1,l_2]_{I,J}$ contains
an element $(\mu,r;\nu,s)$, then we have
\begin{equation}
\rho(I,J)\leq P^{(M)}_{\mu,\nu}[l_1],\quad
\sigma(J)\leq Q^{(N)}_{\mu,\nu}[l_2].
\label{INEQ}
\end{equation}
\end{proposition}
\subsection{Vectors $\rho'$ and $\sigma'$ and upper subsets}\label{UPSEC}

The basic idea in Theorem \ref{recursion theorem} is to change
the rigged partitions with degrees $(M,N-1)$ to those with degrees $(M,N)$.
The parameters $(I,J)$ describes the change of the partitions
{}from $(\mu',\nu')$ given by $m'_\alpha$, $n'_\alpha$
to $(\mu,\nu)$ given by $m_\alpha$, $n_\alpha$:
\begin{equation}
m_\alpha=m'_\alpha+\varepsilon(I)_\alpha,\quad
n_\alpha=n'_\alpha+\varepsilon(J)_\alpha\,\,(1\leq\alpha\leq k).
\label{SHIFT}
\end{equation}
The corresponding change in the riggings is described by the change 
of the upper bounds:
\bea
\Delta r&=&
P^{(M)}_{\mu,\nu}[l_1]-P^{(M)}_{\mu',\nu'}[l_1']
=\kappa(J)-2\kappa(I)+\kappa[l_1'+1,k]-\kappa[l_1+1,k],\label{CAL1}\\
\Delta s&=&
Q^{(N)}_{\mu,\nu}[l_2]-Q^{(N-1)}_{\mu',\nu'}[l_2']
=\kappa(I)-2\kappa(J)+\kappa[1,l_2]+\kappa[l_2'+1,k].\label{CAL2}
\ena
Here $l'_1,l'_2$ are given by (\ref{primes}).
Note that the results are not explicitly dependent on $(\mu',\nu')$
or $(\mu,\nu)$. They are determined only by $I,J,l_1,l_2,l'_1,l'_2$.

The vectors $\rho$ and $\sigma$ give the lower bounds to
the riggings in the lower subsets. We define the upper subsets
by using the shifted lower bounds $\rho'$ and $\sigma'$. Naturally,
the shifts are given by $\Delta r$ and $\Delta s$.

For an $l_1$-admissible pair $(I,J)$
such that $\card(I)=a$ and $\card(J)=b=a+c$, we define the vectors
$\rho'(I,J),\sigma'(I,J)\in\Z^k$ by
\begin{eqnarray}
\rho'(I,J)&=&\rho(I,J)-\Delta r\nonumber\\
&=&\kappa(I)+\kappa[l_1+1,k]-\sum_{i=p+1}^b\kappa(v_i)
-\sum_{i=a+1}^p\kappa(v'_i)-\kappa[l'_1+1,k]\nonumber\\
&=&\kappa(\tilde I)-\kappa[l'_1+1,k],
\label{RHO'}\\
\sigma'(I,J)&=&\sigma(J)-\Delta s\notag\\
&=&\kappa(J)-\kappa(I)-\kappa[l'_2+1,k],\label{SIGMA'}
\end{eqnarray}
where we use $\tilde I$ defined in Section \ref{ADMIS}.

The following is clear {}from (\ref{SIGMA'}) and (\ref{RHO'}).
\begin{lemma}
We have
\begin{equation}
\rho'(I,J),\sigma'(I,J)\geq0\quad\text{and}\quad\rho'(I,J)_k=\sigma'(I,J)_k=0.
\label{LEMPOS}
\end{equation}
\end{lemma}

The following lemma will be used in Section \ref{upper summation section}.
We follow the setting in Lemmas \ref{MIN-1} and \ref{MIN-2}.
\begin{lemma}
We have
\begin{equation}\label{SIGMA-MIN}
\sigma'(I,J_{\rm min})
=\begin{cases}
\Bigl(\kappa[k-b+1,k-c]-\kappa(I)\Bigr)^+
\quad&\text{if $l_1+c\geq k$};\\
\Bigl(\kappa[l_1-a+1,k-c]-\kappa(\tilde I)\Bigr)^+
\quad&\text{if $l_1+c<k$}.
\end{cases}
\end{equation}
\end{lemma}
\begin{proof}
If $l_1+c\geq k$, using (\ref{SIGMA'}) and (\ref{MIN1}) we have
\begin{eqnarray}
\sigma'(I,J_{\rm min})&=&\sum_{i=1}^a\kappa({\rm min}(u_i,k-b+i))-\kappa(I)
\nonumber\\
&=&\Bigl(\kappa[k-b+1,k-c]-\kappa(I)\Bigr)^+.\label{SIG-MIN1}
\end{eqnarray}

If $l_1+c<k$, set $\tilde I=I\sqcup I'$. We have $I'\subset[u_{a+1},k]$ and
$[u_{a+1},k]\backslash\tilde I=[u_{a+1},k]\backslash I'$.
Therefore, using (\ref{SIGMA'}) and (\ref{MIN2}) we have
\begin{eqnarray}
\sigma'(I,J_{\rm min})
&=&\sum_{i=1}^a\kappa({\rm min}(u_i,l_1-a+i))-\kappa(\tilde I)
+\kappa[l_1+1,l_2']
\nonumber\\
&=&\Bigl(\kappa[l_1-a+1,l_2']-\kappa(\tilde I)\Bigr)^+.\label{sigma'2}
\end{eqnarray}
\end{proof}

For $l_1,a,b$ ($0\leq a\leq b$)
and $(I,J)$ such that $\card(I)=a$ and $\card(J)=b$,
we define the subset $R_{m-a,n-b}[l_1]^{I,J}\subset R_{m-a,n-b}$
as follows. If $(I,J)$ is $l_1$-admissible, we set
\begin{eqnarray}\label{UPSET}
&&R_{m-a,n-b}[l_1]^{I,J}=\{
(\mu',r';\nu',s')\in R_{m-a,n-b};\notag\\
&&\quad\rho'(I,J)_\al\leq^{\varepsilon(I)_\al}r'[\al]
\text{\quad\rm and\quad}
\sigma'(I,J)_\al\leq^{\varepsilon(J)_\al}s'[\al]
\text{\quad\rm for all $1\leq\alpha\leq k$}\}
\end{eqnarray}
(see (\ref{UPARROW}) and (\ref{INF})). Otherwise, we set
$R_{m-a,n-b}[l_1]^{I,J}=\emptyset$.

If $l_1+c\geq k$, the set $R_{m-a,n-b}[l_1]^{I,J}$ is independent of $l_1$.
Sometimes we abbreviate $R_{m-a,n-b}[l_1]^{I,J}$  to
$R_{m-a,n-b}^{I,J}$ in this case in avoiding confusion caused by the presence
of $l_1$ in the written formulas.

We call the marking of $\rho',\sigma'$ as before.
The restrictions (\ref{UPSET}) for $r'[\alpha]$ or $s'[\alpha]$
are marked if and only if $\varepsilon(I)_\alpha=-1$ or
$\varepsilon(J)_\alpha=-1$, respectively, and hence they are equalities.

For $M,N-1\geq0$ we define
\begin{equation}
R^{(M,N-1)}_{m-a,n-b}[l_1]^{I,J}=
R_{m-a,n-b}[l_1]^{I,J}\cap R^{(M,N-1)}_{m-a,n-b}[l'_1,l'_2].
\end{equation}

We have
\begin{lemma}\label{NON2}
For $I,J\subset\{1,\dots,k\}$ such that $\card(I)=a,\card(J)=b$,
if $R^{(M,N-1)}_{m-a,n-b}[l_1]^{I,J}$ contains an element $(\mu',r';\nu',s')$,
then we have
\begin{equation}
\rho'(I,J)_\alpha\leq P^{(M)}_{\mu',\nu'}[l'_1]_\alpha,
\quad\sigma'(I,J)_\alpha\leq Q^{(N-1)}_{\mu',\nu'}[l'_2]_\alpha.\label{INEQ2}
\end{equation}
\end{lemma}
\begin{proof}
The proof is completely parallel to Proposition \ref{NON}
(we use (\ref{SIGMA'}) and (\ref{RHO'}))
except that the inequalities (\ref{INEQ2}) for $\alpha=k$ follow
directly {}from (\ref{COC}) and (\ref{LEMPOS}). 
\end{proof}

\subsection{Bijection}\label{BIJECTION}

Define the map 
\be
{\mathfrak m}_{I,J}:R_{m-a,n-b}^{(M,N-1)}[l_1]^{I,J}
\to R_{m,n}^{(M,N)}[l_1,l_2]_{I,J}
\ee
by the formula
\be 
{\mathfrak m}_{I,J}(\mu',r';\nu',s')=(\mu,r,\nu,s),
\ee
where
\begin{equation}
\mu=\mu'+\varepsilon(I),\quad\nu=\nu'+\varepsilon(J),\label{CONN}
\end{equation}
and
the riggings $r,s$ are defined by
\be
\begin{matrix}
r^{(\al)}_i={r'}^{(\al)}_i + (\Delta r)_\al &(1\leq i\leq m'_\al-1);\\
r^{(\al)}_{m'_\al}={r'}^{(\al)}_{m'_\al} + (\Delta r)_\al &
\text{\rm if $\varepsilon(I)_\al=0,1$};\\
r^{(\al)}_{m'_\al+1}=\rho(I,J)_\al &
\text{\rm if $\varepsilon(I)_\al=1$},
\end{matrix}
\ee
and 
\be
\begin{matrix}
s^{(\al)}_i={s'}^{(\al)}_i + (\Delta s)_\al & (1\leq i\leq n'_\al-1);\\
s^{(\al)}_{n'_\al}={s'}^{(\al)}_{n'_\al} + (\Delta s)_\al &
\text{\rm if $\varepsilon(J)_\al=0,1$};\\
s^{(\al)}_{n'_\al+1}=\sigma(J)_\al &
\text{\rm if $\varepsilon(J)_\al=1$}.
\end{matrix}
\ee

We conclude this section by proving
\begin{proposition}
For any $I,J\subset\{i,\dots,k\}$ the map ${\mathfrak m}_{I,J}$ is a bijection.
\end{proposition}
\begin{proof}
It is enough to show the bijectivity of ${\mathfrak m}_{I,J}$
between the subset of $R_{m-a,n-b}^{(M,N-1)}[l_1]^{I,J}$ with a fixed
$\mu',\nu'$ and the subset of $R_{m,n}^{(M,N)}[l_1,l_2]_{I,J}$
with $\mu,\nu$ given by (\ref{CONN}).
Because of Lemma \ref{NON} and Lemma \ref{NON2}, and the definitions
(\ref{CAL1}),(\ref{CAL2}),(\ref{RHO'}) and (\ref{SIGMA'}), these two subsets
are both empty or the inequalities (\ref{INEQ}) and (\ref{INEQ2})
are both valid. In both cases, the bijectivity is clear.
\end{proof}

\section{Decomposition of $R^{(M,N)}_{m,n}[l_1,l_2,l_3]$}
\label{lower summation section}
Fix $k,l_1,l_2$ and $l_3$ as $(\ref{l condition})$.
The aim of this section is to decompose
the set $R^{(M,N)}_{m,n}[l_1,l_2,l_3]$ as
\[
R^{(M,N)}_{m,n}[l_1,l_2,l_3]=
\bigsqcup_{I,J\atop\card(I)\leq l_3,\card(J)\leq l_2}
R^{(M,N)}_{m,n}[l_1,l_2]_{I,J}.
\]
Namely, we decompose the left hand side, in which the riggings $r$ and $s$
are restricted {}from below by the condition (\ref{exception condition}), into
the subsets in the right hand side, in which the riggings are restricted
{}from below separately for each $r[\alpha]$ and $s[\alpha]$ according
to $(I,J)$.

In fact, it is enough to decompose $R_{m,n}[l_1,l_2,l_3]$ as
\[
R_{m,n}[l_1,l_2,l_3]=
\bigsqcup_{I,J\atop\card(I)\leq l_3,\card(J)\leq l_2}
R_{m,n}[l_1,l_2]_{I,J}.
\]

The proof will be carried out in two steps.

The first step is to take the union of the sets
$R_{m,n}[l_1,l_2]_{I,J}$ over $I$ for a fixed $J$.
The is done in Lemma \ref{LEM-LOWSUM}; the union is denoted
by $R_{m,n}[l_1,l_2,l_3]_J$. The idea of the proof is simple.
For a given non-negative integer $t$ the set of integers
$\{i;i\geq t\}$ is the disjoint union of $\{i;i\geq t+1\}$ and $\{i;i=t\}$.
We need more elaborate arguments in the proof.
However, it is done by a successive application of this simple fact.

The second step is to take the union of the sets
$R_{m,n}[l_1,l_2,l_3]_J$ over $J$ and obtain $R_{m,n}[l_1,l_2,l_3]$.
First we carry out this step for $l_3={\rm min}(l_1,l_2)$.
This is actually a special case of the first step.
We obtain $R_{m,n}=R_{m,n}[l_1,l_2,{\rm min}(l_1,l_2)]$ as the union.
Then, we show that the complement in $R_{m,n}$ of the union of
$R_{m,n}[l_1,l_2,l_3]_J$ is equal to the union of its complement in
$R_{m,n}[l_1,l_2,{\rm min}(l_1,l_2)]_J$. This is done by using another simple
fact that the complement $\{i;i\geq 0\}\backslash\{i;i\geq t\}$ is
the union of $\{i;i=s\}$ for $0\leq s\leq t-1$.

\subsection{Union of $R_{m,n}[l_1,l_2]_{I,J}$ over $I$}

We denote $\card(I)=a$ and $\card(J)=b$ as before.
Given $J$ such that $b\leq l_2$, we set $p=p(l_1,J)$ as in Section \ref{ADMIS}.
Define
\begin{eqnarray}
&&I_{\rm max}(J)=\{v_1,\dots,v_{{\rm min}(l_3,p)}\},\quad
\rho_{\rm max}(J)=\rho(I_{\rm max}(J),J),\label{JMAX}\\
&&R_{m,n}[l_1,l_2,l_3]_J
=\{(\mu,r;\nu,s)\in R_{m,n};\nonumber\\
&&r[\alpha]\geq\rho_{\rm max}(J)_\alpha(1\leq\alpha\leq k),
s[\alpha]\geq_{\varepsilon(J)_\alpha}\sigma(J)_\alpha(1\leq\alpha\leq k)\}.
\label{J}
\end{eqnarray}

We also define the subset of $2^{\{1,\dots,k\}}$:
\begin{equation}
T^{(k)}(J;l_1,l_3)=\{I\subset\{1,\dots,k\};
a\leq l_3\,\,\text{and $(I,J)$ is $l_1$-admissible}\}.
\end{equation}
If ${\rm min}(l_3,p)=0$, $T^{(k)}(J;l_1,l_3)=\{\emptyset\}$.
If  ${\rm min}(l_3,p)>0$, we define the structure of colored graph on
$T^{(k)}(J;l_1,l_3)$ as follows.

If $I\in T^{(k)}(J;l_1,l_3)$ and $I\not=I_{\rm max}(J)$,
we draw an outgoing arrow {}from $I$.
We denote the terminal of this arrow by $\xi(I)\in T^{(k)}(J,l_1,l_3)$
and associate the arrow with color $c(I)\in\{1,\dots,{\rm min}(l_3,p)\}$.
The data $\xi(I)$ and $c(I)$ are determined as follows.

Consider $I=\{u_i\}$, $J=\{v_i\}$ and $\{v'_i\}$ as in Section \ref{ADMIS}.
If $u_i=v_i$ for $1\leq i\leq a$, we have $a<{\rm min}(l_3,p)$
since otherwise $I=I_{\rm max}(J)$. We set
$c(I)=a+1\leq{\rm min}(l_3,p)$ and $\xi(I)=I\sqcup\{v'_{c(I)}-1\}$.
Note that $v'_{c(I)}-1\not\in I$ because
$u_a=v_a<v_{c(I)}<v'_{c(I)}$.
If there exists $i$ such that $u_i>v_i$, we set
$c(I)$ to be the minimal integer $i$ satisfying this property, and
$\xi(I)=\Bigl(I\backslash\{u_{c(I)}\}\Bigr)\sqcup\{u_{c(I)}-1\}$.
Note that $u_{c(I)}-1\not\in I$, since otherwise we have a contradiction
$$
u_{c(I)}-1=u_{c(I)-1}=v_{c(I)-1}\leq v_{c(I)}-1<u_{c(I)}-1.
$$

We have
\begin{lemma}\label{LEM-LOWSUM}
\begin{equation}
\bigsqcup_{I\in T^{(k)}(J;l_1,l_3)}R_{m,n}[l_1,l_2]_{I,J}
=R_{m,n}[l_1,l_2,l_3]_J.\label{J1}
\end{equation}
\end{lemma}
\begin{proof}
We use induction on $l_3$. If $l_3=0$, the statement is obvious because the
union (\ref{J1}) is for a single element $I=\emptyset$. We reduce the
proof for $l_1,l_2,l_3,k$ to $l_1-1,l_2-1,l_3-1,k-1$.

Fix $J=\{v_1,\dots,v_{b}\}$ such that $b\leq l_2$,
and denote $R_I=R_{m,n}[l_1,l_2]_{I,J}$.
We take the union of $R_I$ over 
a maximal string $I[i]\in T^{(k)}(J,l_1,l_3)$ ($1\leq i\leq\gamma$)
of color $1$:
$$
I[1]{\buildrel1\over\rightarrow}
I[2]{\buildrel1\over\rightarrow}
\dots{\buildrel1\over\rightarrow}
I[\gamma].
$$
This is maximal in the sense that there is no arrow of color $1$ pointing
to $I[1]$ or {}from $I[\gamma]$. Each arrow of color $1$ belongs to one and only
one maximal string of color $1$.

If $\card(I[\gamma])=1$, $\gamma=v'_1-v_1+1$, $I[1]=\emptyset$ and
$I[i]=\{v'_1-i+1\}$ for $2\leq i\leq\gamma$.
If $a=\card(I[\gamma])>1$, there exists a sequence
$$
u_2<\dots<u_a
$$
such that $\gamma=u_2-v_1$ and $I[i]=\{u_1[i],u_2,\dots,u_a\}$ where
$u_1[i]=u_2-i$. Note that in the case $a=1$, the situation is the same if
we set $u_2=v'_1+1$.

Consider the restriction
$r[\alpha]\geq_{\varepsilon(I[i])_\alpha}\rho(I[i],J)_\alpha$
in $R_{I[i]}$ $(1\leq i\leq\gamma)$. Unless $v_1\leq\alpha\leq u_2-2$,
$\varepsilon(I[i])_\alpha$ and $\rho(I[i],J)_\alpha$ are independent of $i$.

If $v_1\leq\alpha\leq u_2-2$, we have
\begin{eqnarray}
\varepsilon(I[i])_\alpha&=&1\text{ if and only if $i=u_2-\alpha$},\\
\rho(I[i],J)_\alpha&=&
\begin{cases}
\rho(I[1],J)_\alpha&\text{ if $1\leq i\leq u_2-\alpha-1$;}\\
\rho(I[1],J)_\alpha-1&\text{ if $u_2-\alpha\leq i\leq\gamma$.}
\end{cases}
\end{eqnarray}
{}{}From these observations follows that $R_{I[i]}$ ($1\leq i\leq\gamma$)
are disjoint, and the union is characterized by the conditions that
\begin{equation}
r[\alpha]\geq_{\varepsilon(I[\gamma]\backslash\{v_1\})_\alpha}
\rho(I[\gamma],J)_\alpha\,\,(1\leq\alpha\leq k),\quad
s[\alpha]\geq_{\varepsilon(J)_\alpha}
\sigma(J)_\alpha\,\,(1\leq\alpha\leq k).\label{J2}
\end{equation}
Since $\rho(I[\gamma],J)_\alpha=0$ for $1\leq\alpha\leq v_1$,
there is no restriction on $r[\alpha]$ for $1\leq\alpha\leq v_1$.
In particular, there is no restriction for $r[1]$.

Now, we modify the graph. We discard $I(i)$ ($1\leq i\leq\gamma-1$)
{}from $T^{(k)}(J,l_1,l_3)$ and replace
the set $R_{I[\gamma]}$ by the union $R'_{I[\gamma]}$ characterized by
(\ref{J2}).
Carrying out this process for all the maximal strings of color $1$,
we obtain a new graph $T^{(k)}(J;l_1,l_3)'$ and
the sets $R'_I$ ($I\in T^{(k)}(J;l_1,l_3)'$). Observe that
$I=\{u_1,\dots,u_{\card(I)}\}\in T^{(k)}(J;l_1,l_3)'$ satisfies the
restriction $u_1=v_1$ and there is no arrow of color $1$
in $T^{(k)}(J;l_1,l_3)'$.

We see that the graph $T^{(k)}(J;l_1,l_3)'$ is isomorphic to
$T^{(k-1)}(J';l_1-1,l_3-1)$ where $J'=\{v_2-1,\dots,v_b-1\}$.
The isomorphism maps $I$ to $I'=\{u_2-1,\dots,u_a-1\}$ and identifies
the color $c$ in the former with the color $c-1$ in the latter.
We have
$\rho_{l_1,k}(I,J)_\alpha=\rho_{l_1-1,k-1}(I',J')_{\alpha-1}$ for
$2\leq\alpha\leq k$, and
$\varepsilon(I\backslash\{v_1\})_\alpha=1$ if and only if
$\varepsilon(I')_{\alpha-1}=1$.

Therefore, the condition
for $r[\alpha]$ in $R'_I$ is exactly the same as the condition
for $r[alpha-1]$ in the subset $R_{m,n}[l_1-1,l_2-1,l_3-1]_{I',J'}$
at the level $k-1$. Thus we have proved (\ref{J1}).
\end{proof}

\subsection{Union of $R_{m,n}[l_1,l_2,l_3]_J$ over $J$}

Consider the subsets indexed by $J$ such that $\card(J)\leq l_2$ (\ref{J}):
$$
R_{m,n}[l_1,l_2,l_3]_J\subset R_{m,n}.
$$
They are disjoint. In fact, the restrictions on the riggings $s$
given by $\sigma(J)$ and $\varepsilon(J)$ are disjoint
(see Lemma \ref{MINCASE} below).

The goal is to show that the union $R_{m,n}[l_1,l_2,l_3]_J$ over $J$
is equal to $R_{m,n}[l_1,l_2,l_3]$.

If $l_3={\rm min}(l_1,l_2)$, we have
\begin{lemma}\label{MINCASE}
\begin{equation}
\bigsqcup_{J:\card(J)\leq l_2}R_{m,n}[l_1,l_2,{\rm min}(l_1,l_2)]_J
=R_{m,n}.\label{UU}
\end{equation}
\end{lemma}
\begin{proof}
If $l_3={\rm min}(l_1,l_2)$, we have ${\rm min}(l_3,p)=p$ since
$p\leq{\rm min}(l_1,l_2)$. {}From this follows that
$I_{\rm max}(J)=\{v_1,\dots,v_p\}$, and therefore, $\rho(I_{\rm max}(J),J)=0$.
Therefore, there is no restriction on $r[\alpha]$ in
$R_{m,n}[l_1,l_2,{\rm min}(l_1,l_2)]_J$.
We take the union of the riggings $s$ subject to
the restriction on $s[\alpha]$. This is equivalent to the special case of
Lemma \ref{LEM-LOWSUM} where $I,J,l_1,l_2,l_3$ are replaced by
$J,[1,l_2],k,l_2,l_2$, respectively. Therefore, the left hand side of 
(\ref{UU}) is disjoint and the equality holds.
\end{proof}

Set
\begin{eqnarray}
C_1&=&
\bigcup_{1\leq\alpha,\beta\leq k}
\bigcup_{0\leq i+j\leq\tau^{(\alpha,\beta)}[l_1,l_2,l_3]-1}
R^{(\alpha,\beta)}[i,j],\label{EXCEPTIONS}\\
R^{(\alpha,\beta)}[i,j]&=&\{(\mu,r;\nu,s)\in
R_{m,n};r[\alpha]=i,s[\beta]=j\}.
\label{DEF-EXCEPTIONS}
\end{eqnarray}
It is easy to see that
$$
C_1=R_{m,n}\backslash R_{m,n}[l_1,l_2,l_3].
$$

Set
\begin{eqnarray}
C_2&=&
R_{m,n}\backslash U,\label{COM}\\
U&=&\bigsqcup_{J:\card(J)\leq l_2}R_{m,n}[l_1,l_2,l_3]_J.\label{U}
\end{eqnarray}

Lemma \ref{MINCASE} enables us to represent $C_2$,
which is by definition the complement of union,
as the union of complements. Namely, we have
\begin{eqnarray}
C_2&=&\bigsqcup_{J:\card(J)\leq l_2} R^c_J,\label{JSUM}\\
R^c_J&=&
R_{m,n}[l_1,l_2,{\rm min}(l_1,l_2)]_J\backslash
R_{m,n}[l_1,l_2,l_3]_J.
\label{DEF-JSUM}
\end{eqnarray}

The goal is to show that
\[
C_1=C_2.
\]

First we assume that $l_3=0$. In this case, we have $I_{\rm max}(J)=\emptyset$.
We prove that (\ref{JSUM}) is equal to (\ref{EXCEPTIONS}).

We call $K\subset[1,k]$ of the first kind
if for some $\beta(K),b(K)\in\{1,\dots,k\}$ it is
of the form
\begin{equation}
K=[\beta(K)-b(K)+1,\beta(K)].\label{OFK}
\end{equation}

We will modify (\ref{JSUM}) and obtain another representation of the form
\begin{eqnarray}
C_3&=&
\bigcup_{K\text{\rm : of the first kind and $\card(K)\leq l_2$}}
R^{'c}_K,\label{REWRITE}\\
R^{'c}_K&=&
\bigcup_{1\leq\alpha\leq k\atop0\leq i\leq\rho(\emptyset,K)_\alpha-1}
R^{'c}_K[\alpha,i],\label{DEF1}\\
R^{'c}_K[\alpha,i]&=&
\{(\mu,r;\nu,s)\in
R_{m,n};\nonumber\\
&&\quad r[\alpha]=i,s[\beta(K)]=\sigma(K)_{\beta(K)}\}.\label{DEF2}
\end{eqnarray}

We start {}from a lemma on some property of the restriction (\ref{J})
on the riggings $s$ given by $\sigma(J)=\kappa[1,l_2]-\kappa(J)$
and $\varepsilon(J)$.

For $J\subset[1,k]$ such that $\card(J)\leq l_2$ we set
\begin{equation}
S_J=\{s=(s_1,\dots,s_k)\in{\Z}_{\geq0}^k;
s_\alpha\geq_{\varepsilon(J)_\alpha}\sigma(J)_\alpha\,\,(1\leq\alpha\leq k)\},
\end{equation}
and for $K=[\beta-b+1,\beta]\subset[1,k]$ such that $b\leq l_2$
\begin{equation}
S'_K=\{s\in{\Z}_{\geq0}^k;s_\beta=\sigma(K)_\beta\}.
\end{equation}
As we have already mentioned in the proof of Lemma \ref{MINCASE},
the subsets $S_J$ are disjoint.

\begin{lemma}\label{INC}
We have the inclusion
\begin{equation}
S'_K\subset\cup_{J\geq K}S_{J}.
\end{equation}
\end{lemma}
\begin{proof}
We will prove this by induction on $K$ with respect to the ordering
defined in Section \ref{ADMIS}.
We see that the statement is true for the maximal element $K=[1,l_2]$.
In fact, if $K=[1,l_2]$ the statement $S'_K=S_K$ follows {}from $\sigma(K)=0$
and $\varepsilon([1,l_2])_\alpha=1$ if and only if $\alpha=l_2$.
This is the base of the induction.

Now assume that the statement is true for all $K'$ of the first kind
such that $K'> K$. We will show that there exists a subset ${\overline S}_K$
satisfying
\begin{eqnarray}
{\overline S}_K&\subset&\bigsqcup_{J\geq K}S_{J},\label{S1}\\
S'_K\backslash{\overline S}_K&\subset&
\bigcup_{K'>K\atop K'\text{\rm : of the first kind}}S'_{K'},\label{S2}
\end{eqnarray}
This will close the induction steps.

We fix $K=[\beta-b+1,\beta]$ and define
\begin{equation}
{\overline S}_K
=\bigsqcup_{\beta\in J\geq K\atop\card(J\cap[1,\beta-1])=b-1}S_J.
\end{equation}
Namely, we take the disjoint union over
$J=\{v_1,\dots,v_{b-1},\beta,v_{b+1},\dots,v_{b'}\}$
such that $1\leq v_1<\dots<v_{b-1}<\beta<v_{b+1}<\dots<v_{b'}$ with
$b'\leq l_2$. Note that the element $\beta$ is fixed, $v_1,\dots,v_{b-1}$
move around the interval $[1,\beta-1]$ and new elements
$v_{b+1},\dots,v_{b'}$ are added in the interval $[\beta+1,k]$.
We have (\ref{S1}) obviously.

By the same argument as in the proof of Lemma \ref{LEM-LOWSUM} we obtain
\begin{equation}
{\overline S}_K
=\{s\in{\Z}_{\geq0}^k;
s_\alpha\geq\sigma(K_{\rm max})_\alpha\,\,(1\leq\alpha\leq k)\},
\end{equation}
where $K_{\rm max}=[1,b-1]\sqcup[\beta,\beta+l_2-b]$.

We have the following values of $\sigma(K)_\alpha$
and $\sigma(K_{\rm max})_\alpha$.

If $\beta\geq l_2$, then
\begin{eqnarray}
\sigma(K)_\alpha&=&
\begin{cases}
\alpha&(1\leq\alpha\leq l_2);\\
l_2&(l_2\leq\alpha\leq\beta-b);\\
l_2+\beta-\alpha-b&(\beta-b\leq\alpha\leq\beta);\\
l_2-b&(\beta\leq\alpha\leq k),
\end{cases}\\
\sigma(K_{\rm max})_\alpha&=&
\begin{cases}
{\rm max}(0,\alpha-b+1)&(1\leq\alpha\leq l_2);\\
l_2-b+1&(l_2\leq\alpha\leq\beta-1);\\
{\rm max}(0,l_2+\beta-\alpha-b)&(\beta-1\leq\alpha\leq k).\\
\end{cases}
\end{eqnarray}

If $\beta\leq l_2$, then
\begin{eqnarray}
\sigma(K)_\alpha&=&
\begin{cases}
\alpha&(1\leq\alpha\leq\beta-b);\\
\beta-b&(\beta-b\leq\alpha\leq\beta);\\
\alpha-b&(\beta\leq\alpha\leq l_2);\\
l_2-b&(l_2\leq\alpha\leq k),
\end{cases}\\
\sigma(K_{\rm max})_\alpha&=&
\begin{cases}
{\rm max}(0,\alpha-b+1)&(1\leq\alpha\leq\beta-1);\\
\beta-b&(\beta-1\leq\alpha\leq l_2);\\
{\rm max}(0,l_2+\beta-\alpha-b)&(l_2\leq\alpha\leq k).\\
\end{cases}
\end{eqnarray}

Now we will prove  (\ref{S2}). We have
$\sigma(K)_\beta=\sigma(K_{\rm max})_\beta$. Therefore,
\begin{equation}
S'_K\backslash{\overline S}_K
=\bigcup_{\alpha\not=\beta\atop0\leq i\leq\sigma(K_{\rm max})_\alpha-1}
\{s\in{\Z}_{\geq0}^k;s_\alpha=i\}\cap S'_K.\label{30}
\end{equation}
We take $K'$ in (\ref{S2}) to be $K_{\alpha,b'}=[\alpha-b'+1,\alpha]$. We have
$K_{\alpha,b'}\subset[1,k]$ and $K_{\alpha,b'}\geq K$ if and only if
\begin{equation}
{\rm max}(b,b+\alpha-\beta)\leq b'\leq{\rm min}(l_2,\alpha).\label{REGION}
\end{equation}
By case checking one can prove that the set of integers
consisting of the values
of $\sigma(K_{\alpha,b'})_\alpha={\rm min}(l_2,\alpha)-b'$ where $b'$ runs over
(\ref{REGION}), contains $[0,\sigma(K_{\rm max})_\alpha-1]$
appearing in (\ref{30}). For example,
if $1\leq\alpha\leq l_2\leq\beta$, we have
$\sigma(K_{\rm max})_\alpha={\rm max}(0,\alpha-b+1)$ and
$\sigma(K_{\alpha,b'})_\alpha=\alpha-b'$. Therefore, we obtain
$$
\cup_{0\leq i\leq\sigma(K_{\rm max})_\alpha-1}
\{s\in{\Z}_{\geq0}^k;s_\alpha=i\}\cap S'_K
\subset\cup_{b\leq b'\leq\alpha}S'_{K_{\alpha,b'}}.
$$
Other cases are similar.
\end{proof}

Now we prove
\begin{lemma}\label{C2C3}
Assume that $l_3=0$. We have
\[
C_2=C_3.
\]
\end{lemma}
\begin{proof}
If $l_3=0$ we have
\begin{eqnarray}
R^c_J&=&\bigcup_{1\leq\alpha\leq k\atop0\leq i\leq\rho(\emptyset,J)_\alpha-1}
R^c_J[\alpha,i],\nonumber\\
R^c_J[\alpha,i]
&=&\{(\mu,r;\nu,s)\in R_{m,n};\nonumber\\
&&r[\alpha]=i,s[\beta]\geq_{\varepsilon(J)_\beta}\sigma(J)_\beta
\,\,(1\leq\beta\leq k)\}.\label{RCJ}
\end{eqnarray}

First we show that if $K$ is an interval of the first kind (\ref{OFK}) we have
\begin{equation}
R^{'c}_K\subset\cup_{J\geq K}R^c_{J}.
\label{FIRST}
\end{equation}
>{}From this follows that $C_3\subset C_2$.

If $J'>J$ then $\rho(\emptyset,J')\geq\rho(\emptyset,J)$.
Therefore, in order to show (\ref{FIRST})
one can forget the restriction on $r$. Then, it follows {}from Lemma \ref{INC}.

To finish the proof, we show that $C_2\subset C_3$.
Consider $J=J^{(1)}\sqcup\dots\sqcup J^{(h)}$ where
$$
J^{(j)}=[\beta^{(j)}-b^{(j)}+1,\beta^{(j)}]
$$
and $\beta^{(j)}<\beta^{(j+1)}-b^{(j+1)}$.
Set
$$
K^{(j)}=[\beta^{(j)}-(b(1)+\dots+b(j))+1,\beta^{(j)}]\quad(1\leq j\leq h).
$$
We will show that
\begin{equation}
R^c_J\subset\cup_{j=1}^hR^{'c}_{K^{(j)}}.\label{STATEMENT}
\end{equation}
Note that $\sigma(J)_{\beta^{(j)}}=\sigma(K^{(j)})_{\beta^{(j)}}$
and $\varepsilon(J)_{\beta^{(j)}}=1$. Therefore, the condition for the
riggings $s$ in $R^c_J[\alpha,i]$ is stronger than
$R^{'c}_{K^{(j)}}[\alpha,i]$.
Namely, we have $R^c_J[\alpha,i]\subset R^{'c}_{K^{(j)}}[\alpha,i]$.

For any $\alpha$ we can find $j$ such that
$$
\beta^{(j)}-b^{(j)}+1\leq\alpha\leq\beta^{(j+1)}-b^{(j+1)}.
$$
Then we have
$$
\rho(\emptyset,J)_\alpha=\rho(\emptyset,K^{(j)})_\alpha.
$$
The statement (\ref{STATEMENT}) follows {}from this.
\end{proof}

Next we show
\begin{lemma}\label{C1C3}
Assume that $l_3=0$. We have
\[
C_1=C_3.
\]
\end{lemma}
\begin{proof}
We show that (\ref{REWRITE}) is equal to (\ref{EXCEPTIONS}) by case checking
for each case of the ordering of $\alpha,\beta,l_1,l_2$. There are 24 cases.
Here we give the details for the case
$l_2\leq l_1<\alpha\leq\beta$. Other cases are similar.

If $l_2\leq l_1<\alpha\leq\beta$ the intervals $K$ which appear in
(\ref{REWRITE}) satisfying $\beta(K)=\beta$ are of the form $J=[\gamma,\beta]$
where $\beta-l_2+1\leq\gamma\leq\alpha$. For such $K$ we have
$$
\rho(\emptyset,K)_\alpha=l_1-\gamma+1,\quad
\sigma(K)_\beta=l_2-(\beta-\gamma+1).
$$
Therefore, the pair of integers
$(\rho,\sigma)=(\rho(\emptyset,K)_\alpha-1,\sigma(K)_\beta)$
runs over the set
$\{(\rho,\sigma);\rho,\sigma\geq0,\rho+\sigma=l_1+l_2-\beta-1\}$.
On the other hand we have
$$
\tau^{(\alpha,\beta)}[l_1,l_2,0]={\rm min}
(\alpha,\beta,l_1,l_2,l_1+\beta-\alpha,,l_2+\alpha-\beta,
l_1+l_2-\alpha,l_1+l_2-\beta)=l_1+l_2-\beta.
$$
This completes the proof.
\end{proof}

Finally, we have
\begin{lemma}
\begin{equation}\label{JUNION}
\bigsqcup_{J:\card(J)\leq l_2}R_{m,n}[l_1,l_2,l_3]_J
=R_{m,n}[l_1,l_2,l_3].
\end{equation}
\end{lemma}
\begin{proof}
By Lemmas \ref{C2C3} and \ref{C1C3} we have shown (\ref{JUNION}) for $l_3=0$.
Let us reduce the proof to the case $l_3=0$.
Suppose that $l_3>0$. Then, we have $l_1,l_2>0$.
We will reduce this case to the case where
$l_1,l_2,l_3$ and $k$ replaced by $l_1-1,l_2-1,l_3-1$ and $k-1$, respectively.

Note that the union is taken over $J$ such that $\card(J)\leq l_2$, i.e.,
$J\in T^{(k)}([1,l_2];k,l_2)$. Therefore, we
refer to the structure of colored graph in this set.

Recall the definition of $\rho_{\rm max}(J)$ given by (\ref{JMAX}).
If $J$ varies on a maximal string of color $1$, then
only $v_1$ changes. However, we see that
the $\rho_{\rm max}(J)$ is independent of $v_1$ because for $I=I_{\rm max}(J)$
we have $u_1=v_1$. It follows that
the vector $\rho_{\rm max}(J)$ is constant on the maximal string.
Therefore, we can take the union over $J$ on
maximal strings of color $1$ only on the riggings $s$ forgetting $r$.

Taking unions over all of the maximal strings of color $1$,
we can rewrite the left hand side of (\ref{JUNION})
as the union of the resulting subsets over such $J$
that satisfies $1\in J$, i.e., of the form $J=\{1,v_2,\dots,v_b\}$.
The subgraph of $T^{(k)}([1,l_2];k,l_2)$ consisting of such $J$ is
isomorphic to $T^{(k-1)}([1,l_2-1];k-1,l_2-1)$ by mapping
$J$ to $J'=\{v_2-1,\dots,v_b-1\}$ and identifying the color $c$
in the former with the color $c-1$ in the latter.

We have
$$
\rho_{l_1,k}(I_{\rm max}(J),J)_\alpha
=\rho_{l_1-1,k-1}(I_{\rm max}(J'),J')_{\alpha-1}\quad(2\leq\alpha\leq k).
$$
Therefore, we have
\begin{eqnarray*}
\sigma_{l_2,k}(J)_\alpha&=&(\kappa([1,l_2])-\kappa(J))_\alpha\\
&=&(\kappa([1,l_2-1])-\kappa(J'))_{\alpha-1}\\
&=&\sigma_{l_2-1,k-1}(J')_{\alpha-1}.
\end{eqnarray*}
Note also that
\[
\tau^{(\alpha,\beta)}[l_1,l_2,l_3]=
\tau^{(\alpha-1,\beta-1)}[l_1-1,l_2-1,l_3-1].
\]
Thus, we have reduced the case $l_1$, $l_2$, $l_3$, $k$ to
$l_1-1,l_2-1,l_3-1,k-1$.
\end{proof}

In conclusion, we have
\begin{proposition}\label{CONCL}
\begin{equation}
R_{m,n}[l_1,l_2,l_3]=
\bigsqcup_{\card(I)\leq l_3,\card(J)\leq l_2}
R_{m,n}[l_1,l_2]_{I,J}.
\label{LOWSUM-2}
\end{equation}
\end{proposition}
\begin{proposition}
\begin{equation}
R^{(M,N)}_{m,n}[l_1,l_2,l_3]=
\bigsqcup_{\card(I)\leq l_3,\card(J)\leq l_2}
R^{(M,N)}_{m,n}[l_1,l_2]_{I,J}.
\end{equation}
\end{proposition}

\section{Decomposition of $R^{(M,N-1)}_{m-a,n-b}[l'_1,l'_2,l'_3]$}
\label{upper summation section}
Fix $k,l_1,a,b$ such that $0\leq a\leq l_1\leq k$, $a\leq b=a+c\leq k$,
and define $l'_1,l'_2,l'_3$ by (\ref{primes}). We denote
$m'=m-a$ and $n'=n-b$ in this section.

The aim of this section is to decompose the set
$R_{m',n'}^{(M,N-1)}[l_1',l_2',l_3']$ as
\begin{equation}
R_{m',n'}^{(M,N-1)}[l_1',l_2',l_3']=
\bigsqcup_{\card (I)=a,\card(J)=b}
R_{m',n'}^{(M,N-1)}[l_1]^{I,J}.
\end{equation}
Again, it is enough to decompose $R_{m',n'}[l_1',l_2',l_3']$ as
\begin{equation}
R_{m',n'}[l_1',l_2',l_3']=
\bigsqcup_{\card (I)=a,\card(J)=b}R_{m',n'}[l_1]^{I,J}.
\label{UPDEC}
\end{equation}

It is useful to note that if $\alpha=k$ or $\beta=k$ then
$\tau^{(\alpha,\beta)}[l'_1,l'_2,l'_3]=0$. Also, because of (\ref{LEMPOS})
there are no restrictions on $r'[k]$ nor $s'[k]$ in the definition
(\ref{UPSET}) of the upper subsets $R_{m',n'}[l_1]^{I,J}$.
Therefore, we can restrict our discussion on $k$ vectors in this section
to the interval $1\leq\alpha\leq k-1$.

The proof is divided into two cases: $l_1+c\geq k$ and $l_1+c<k$.
\subsection{Case $l_1+c\geq k$} 

In this case, we have (see \Ref{primes})
\be
l_1'=k-a,\quad l_2'=k-c\quad l'_3=k-b.
\label{THATCASE}
\ee
and
\begin{align*}
&\sigma'(I,J)=\kappa(J)-\kappa(I)-\kappa[k-c+1,k],\\
&\rho'(I,J)=\kappa(I)-\kappa[k-a+1,k].
\end{align*}
Note, in particular, that the $l_1$-dependence disappears.
We write $R_{m',n'}^{I,J}$ for $R_{m',n'}[l_1]^{I,J}$.

First, we fix the subset $I=\{u_1,\dots,u_a\}$ and take the union over
$J=\{v_1,\dots,v_b\}$. This is similar to Lemma \ref{LEM-LOWSUM}
We use (\ref{MIN1}) for $J_{\rm min}$, (\ref{RHO'}) with $\tilde I=I$
for $\rho'(I)=\rho'(I,J_{\rm min})$ and (\ref{SIGMA-MIN})
for $\sigma'(I)=\sigma'(I,J_{\rm min})$. They are all independent of $l_1$.

Set
\begin{eqnarray}
R_{m',n'}^I
&=&\{(\mu',r';\nu',s')\in R_{m',n'};\nonumber\\
&&r'[\alpha]\geq^{\varepsilon(I)_\alpha}\rho'(I)_\alpha
(1\leq\alpha\leq k-1),
s'[\alpha]\geq\sigma'(I)_\alpha
(1\leq\alpha\leq k-1)\}.\nonumber\\\label{I}
\end{eqnarray}

\begin{lemma}\label{fixed-I-summation} 
For $I=\{u_1,\dots,u_a\}$, we have
\be
\bigsqcup_{ J={\{v_1,\dots,v_b\}\subset\{1,\dots,k\}} 
\atop {v_1\leq u_1,\dots, v_a\leq u_a} }
R_{m',n'}^{I,J}=R_{m',n'}^I.
\label{LOW-J-SUM}
\ee
\end{lemma}
\begin{proof}
The proof of Lemma \ref{fixed-I-summation} is parallel to that of Lemma 
\ref{LEM-LOWSUM}.

In Lemma \ref{fixed-I-summation} the restriction on $r'$ in $R_{m',n'}^{I,J}$
is independent of $J$ and the restriction on $s'$ is of the form
\begin{eqnarray*}
&s'[\alpha]\geq(\kappa(J)+A)_\alpha\quad
&\text{if $\varepsilon(J)_\alpha\not=-1$},\\
&s'[\alpha]=(\kappa(J)+A)_\alpha\quad
&\text{if $\varepsilon(J)_\alpha=-1$}.
\end{eqnarray*}
Here $A$ is a $k$-vector independent of $J$.

In Lemma \ref{LEM-LOWSUM} the restriction on $s$ in $R_{m,n}[l_1,l_2]_{I,J}$
is independent of $I$ and
the restriction on $r$ is of the form
\begin{eqnarray*}
&r[\alpha]\geq(-\kappa(I)+B)_\alpha\quad
&\text{if $\varepsilon(I)_\alpha\not=1$},\\
&r[\alpha]=(-\kappa(I)+B)_\alpha\quad
&\text{if $\varepsilon(I)_\alpha=1$}.
\end{eqnarray*}
Here $B$ is a $k$-vector independent of $I$.

In Lemma \ref{fixed-I-summation} the union is taken over $J$ such that
\[
J_{\rm min}\leq J\leq[1,b],
\]
where $J_{\rm min}$ is given by (\ref{MIN1}), and 
in Lemma \ref{LEM-LOWSUM} the union is taken over $I$ such that
\[
\emptyset\leq I\leq I_{\rm max},
\]
where $I_{\rm max}$ is given by (\ref{JMAX}).

Recall that in the proof of Lemma \ref{LEM-LOWSUM} we take the union over the
maximal strings of color $1$ as the first inductive step.
Similarly, in the setting of Lemma \ref{fixed-I-summation} we take the union
over strings $J[i]$ ($1\leq i\leq\gamma$) of the form
$J[i]=\{v_1,\dots,b_{b-1}[i]\}$ ($v_1<\dots<v_{b-1}$)
and $v_b[i]=v_{b-1}+i$. Here
\[
\gamma=
\begin{cases}
k-v_{b-1}\quad&\text{if $b>a$};\\
u_a-v_{b-1}\quad&\text{if $b=a$}.
\end{cases}
\]

The difference between two cases is that $\card(J)=b$ is fixed in Lemma
\ref{fixed-I-summation}, while $\card(I)$ varies in Lemma \ref{LEM-LOWSUM}.
However, if we consider $\overline I=I\sqcup\{v'_p,\dots,v'_{\card(I)+1}\}$
instead of $I$, $\card(\overline I)=p$ is fixed and two cases are
completely parallel.

Therefore, the union is obtained by substituting $J$ by $J_{\rm min}$
and make the restriction on $s'$ unmarked.
\end{proof}

Now we translate the formula (to be proved)
\begin{equation}\label{TOBE}
\bigsqcup_{\card(I)=a}R_{m',n'}^I=R_{m',n'}[l_1',l_2',l_3']
\end{equation}
into the formula
\begin{equation}
R_{m',n'}[l_1,l_2,0]=\bigsqcup_{\card(J)\leq l_2}
R_{m',n'}[l_1,l_2]_{\emptyset,J},
\label{TRANS}
\end{equation}
which is the special case of (\ref{LOWSUM-2}) with $l_3=0$. We use the case
\[
\text{$l_1=c$ and $l_2=a$}
\]
by the following reason.

In $R_{m',n'}^I$ we have
\begin{eqnarray}
\rho'(I)&=&\kappa(I)-\kappa[k-a+1,k]
\,\,=\,\,\kappa(I)-\kappa[l'_1+1,k],\label{COMP1}\\
\sigma'(I)&=&\Bigl(\kappa[k-b+1,k-c]-\kappa(I)\Bigr)^+
\,\,=\,\,\Bigl(\kappa[l'_3+1,l'_2]-\kappa(I)\Bigr)^+.\label{COMP2}
\end{eqnarray}
On the other hand, in (\ref{TRANS}) we have
\begin{eqnarray*}
\rho(\emptyset,J)
&=&\Bigl(\kappa(J)-\kappa[c+1,c+\card(J)\Bigr)^+,\\
\sigma(J)&=&\kappa[1,a]-\kappa(J).
\end{eqnarray*}
Note, in particular, that $\rho(\emptyset,J)_k=\sigma(J)_k=0$.

We define an involution of the set $\{1,\dots,k\}$ by
\[
i^\dagger=k+1-i.
\]
Using
\begin{equation}
\kappa(i)_\alpha+\kappa(i^\dagger)_{k-\alpha}=1,
\end{equation}
we obtain
\begin{eqnarray}
\rho'(I)_\alpha=\sigma(I^\dagger)_{k-\alpha},\label{RS1}\\
\sigma'(I)_\alpha=\rho(\emptyset,I^\dagger)_{k-\alpha}.\label{RS2}
\end{eqnarray}

In this way, we can translate (\ref{TOBE}) into (\ref{TRANS}) except that
the union is taken over $I$ with the fixed size $\card(I)=a$ in (\ref{TOBE})
while $J$ in (\ref{TRANS}) is only restricted by $\card(J)\leq a$.

Therefore, we need some modification. We will take the union in (\ref{TRANS})
partially so that only $J$ of size $a$ remain.

Given $J=\{v_1,\dots,v_{a'}\}$ such that $a'\leq a$ we define
the closure of $J$ by
\begin{equation}
\overline J=J\sqcup\{w_1,\dots,w_{a-a'}\}
\end{equation}
where $w_1<\dots<w_{a-a'}$ are chosen to be the maximal $a-a'$
elements in $[1,k]\backslash J$.

For $K$ such that $\card(K)=a$ we set
\begin{eqnarray}
\bigl(R_{m',n'}\bigr)_K&=&\{(\mu,r;\nu,s)\in R_{m',n'};\nonumber\\
&&r[\alpha]\geq\rho(\emptyset,K)_\alpha\quad(1\leq\alpha\leq k-1),\nonumber\\
&&s[\alpha]\geq_{\varepsilon(K)}\sigma(K)_\alpha\quad(1\leq\alpha\leq k-1)\}.
\end{eqnarray}
Note that we impose no restrictions at $\alpha=k$.

We have
\begin{lemma}\label{REDUCARG}
\begin{equation}
\bigl(R_{m',n'}\bigr)_K
=\bigsqcup_{J:\overline J=K}R_{m',n'}[c,a]_{\emptyset,J}.
\end{equation}
\end{lemma}
\begin{proof}
Suppose that $\card(J)=a'\leq a$ and $\overline J=K=K_0\sqcup[k-d+1,k]$
where $k-d\not\in K_0$. Then, we have $K_0\subset J$,
and $K_0$ and $d$ are uniquely determined {}from $K$.
We have
\begin{eqnarray*}
\rho(\emptyset,J)&=&(\kappa(J)-\kappa[c+1,c+a'])^+\\
&=&(\kappa(K_0)-\kappa[c+1,c+a-d]
+\kappa(J\backslash K_0)-\kappa[c+a-d+1,c+a'])^+.
\end{eqnarray*}
Note that $\card(K_0)=a-d$ and $\card(J\backslash K_0)=a'+d-a$.
Therefore, we have $(\kappa(K_0)-\kappa[c+1,c+a-d])_\alpha\leq0$
for $\alpha\geq c+a-d$. Since $J\backslash K_0\subset[k-d+1,k]$, we have
\[
\kappa(J\backslash K_0)\leq\kappa[k-d+1,k+a'-a]\leq \kappa[c+a-d+1,c+a'].
\]
>{}From these observations follows that
\begin{eqnarray*}
\rho(\emptyset,J)&=&(\kappa(K_0)-\kappa[c+1,c+a-d])^+\\
&=&\rho(\emptyset,K).
\end{eqnarray*}
The set $J$ satisfies $\card(J)\leq a$ and $\overline J=K$
if and only if $J=K_0\sqcup J'$ where $J'\subset [k-d+1,k]$.

By a similar argument as the proof of Lemma \ref{LEM-LOWSUM},
taking the union of $R_{m',n'}[c,a]_{\emptyset,J}$ over $J'$,
we obtain $(R_{m',n'})_K$.
\end{proof}

Now we are ready to finish the proof of (\ref{UPDEC}) for $l_1+c\geq k$.
\begin{lemma}
\begin{equation}
\bigsqcup_{\card(I)=a}R^I_{m',n'}=R_{m',n'}[l'_1,l'_2,l'_3]
\label{N-1SUM}
\end{equation}
\end{lemma}
\begin{proof}
Let us observe that there is a correspondence between $R^I_{m',n'}$
and $(R_{m',n'})_K$. We have (\ref{RS1}) and (\ref{RS2}).
Moreover, because of
\begin{eqnarray}
\varepsilon(i)_\alpha+\varepsilon(k+1-i)_{k-\alpha}=0,\label{EQEPS}
\end{eqnarray}
$\varepsilon(I^\dagger)_{k-\alpha}=1$ 
if and only if $\varepsilon(I)_\alpha=-1$.

Finally, note that 
\be
\tau^{(k-\beta,k-\al)}[k-a,k-c,k-a-c]=\tau^{(\al,\beta)}[c,a,0].
\ee
In this way, (\ref{N-1SUM}) follows {}from Proposition \ref{CONCL}.
\end{proof}
\subsection{Case $l_1+c<k$} 

In this case, we have (see \Ref{primes})
\be
l_1'=l_1+c-a,\quad l_2'=k-c\quad l'_3=l_1-a.
\label{THISCASE}
\ee
We use (\ref{MIN2}) for $J_{\min}$, (\ref{RHO'}) for 
$\rho'(I)=\rho'(I,J_{\rm min})$ and (\ref{SIGMA'}) for
$\sigma'(I)=\sigma'(I,J_{\rm min})$. Namely, we have
\bea
\rho'(\tilde I)
&=&\kappa(\tilde I)-\kappa[l_1'+1,k],\label{RHO'2}\\
\sigma'(\tilde I)
&=&(\kappa[l'_3+1,l_2']-\kappa(\tilde I))^+.\label{SIGMA'2}
\ena

Because $I'\subset[l_1+1,k]$, we have
\begin{equation}
(\kappa[l_1-a+1,l_2']-\kappa(\tilde I))^+=(\kappa[l_1-a+1,l_1]-\kappa(I))^+
+\kappa[l_1+1,l'_2]-\kappa(I').\label{POSITIVITY}
\end{equation}

In the following lemma, when $l_1+c<k$, we define $R^{\tilde I}_{m',n'}$
differently {}from $R^I_{m',n'}$ defined in (\ref{I}) when $l_1+c\geq k$.

\begin{lemma}\label{J-summation}  Suppose that $\tilde I$ satisfies
the condition $(\ref{TIL-I})$. Set
\begin{eqnarray}
\tilde R_{m',n'}^{\tilde I}
&=&\{(\mu',r';\nu',s')\in
R_{m',n'};\nonumber\\
&&r'[\alpha]\geq^{\varepsilon(I)_\alpha}\rho'(\tilde I)_\alpha
(1\leq\alpha\leq k-1),
s'[\alpha]\geq_{\varepsilon(I')_\alpha}\sigma'(\tilde I)_\alpha
(1\leq\alpha\leq k-1)\}.
\label{64}
\end{eqnarray}

We have
\be
\bigsqcup_{\tilde J}
R_{m',n'}[l_1]^{{\mathfrak c}(\tilde I,\tilde J)}=
\tilde R_{m',n'}^{\tilde I}
\ee
where the summation in the LHS is taken over all $\tilde J$ satisfying
(\ref{TIL-J}) (see (\ref{IJ}) for $(I,J)={\mathfrak c}(\tilde I,\tilde J)$).
\end{lemma}
\begin{proof}
The proof of this lemma is parallel to the proof of Lemma
\ref{fixed-I-summation}. Note that $\rho'(\tilde I)$ (see (\ref{RHO'})) does
not depend on $\tilde J$.  The summation extends to all $\tilde J$ satisfying
$$
[1,a+u_{a+1}-l_1-1]\geq\tilde J\geq\tilde J_{\rm min}
$$
where
\be
\tilde J_{\rm min}=\{{\rm min}(u_i,l_1-a+i)\}_{1\leq i\leq a}
\sqcup[l_1+1,u_{a+1}-1].
\ee
The restriction on $s[\alpha]$ in
$R_{m',n'}[l_1]^{{\mathfrak c}(\tilde I,\tilde J)}$ is given by 
$\sigma'(I,J)_\alpha$. It is marked if and only if $\varepsilon(J)_\alpha=-1$.
The restriction on $s[\alpha]$ in the union over $\tilde J$ is given by
$\sigma'(\tilde I)_\alpha$. The marking will change as follows.
Because of $J\cap[1,u_{a+1}-1]=\tilde J$ the restriction on $s[\alpha]$ is
unmarked if $\alpha\in[1,u_{a+1}-1]$.
In the interval $\alpha\in[u_{a+1},k]$, the marking is unchanged because
$J'=J\cap[u_{a+1},k]$ is fixed.

Decompose the interval $[u_{a+1},k]=I'\sqcup J'$
into subintervals $I'_1,\dots,I'_h$
constituting $I'$ and $J'_1,\dots,J'_h$ constituting $J'$
in such a way that ${\rm max}(I'_i)+1={\rm min}(J'_i)$
and ${\rm max}(J'_i)+1={\rm min}(I'_{i+1})$. Note that $u_{a+1}\in I'_1$
and $J'_h$ may be empty. The restriction on $s[\alpha]$ is marked
if and only if $\alpha={\rm min}(J'_i)-1$ for some $1\leq i\leq h$. This is
equivalent to say that it is marked if and only if
$\alpha={\rm max}(I'_i)$ for some $1\leq i\leq h$ except for $\alpha=k$.
Therefore, we have the marking given in (\ref{64}).
\end{proof}

We remake the union over $\tilde I$ as follows.
\begin{lemma}
Set 
\begin{eqnarray}
&&\overline R_{m'.n'}^{\tilde I}=
\{(\mu',r';\nu',s')\in R_{m',n'};
\nonumber\\
&&r'[\alpha]\geq^{\varepsilon(\tilde I)_\alpha}\rho'(\tilde I)_\alpha
\quad(1\leq\alpha\leq k-1;\alpha\not=l_1\text{ if $u_{a+1}=l_1+1$}),\nonumber\\
&&r'[l_1]\geq\rho'(\tilde I)_{l_1}\text{ if $u_{a+1}=l_1+1$},
\nonumber\\
&&s'[\alpha]\geq\sigma'(\tilde I)_\alpha\quad(1\leq\alpha\leq k-1)\}.
\label{TILDEI}
\end{eqnarray}
We have
\begin{equation}
\bigsqcup_{\tilde I}\tilde R_{m',n'}^{\tilde I}
=\bigsqcup_{\tilde I}\overline R_{m',n'}^{\tilde I}.
\label{REVERSE}
\end{equation}
\end{lemma}
\begin{proof}
Fix $\tilde I$. First, we can rewrite each subset $\tilde R_{m',n'}^{\tilde I}$
in terms of subsets which have the same restriction on $r'$
as $\tilde R_{m',n'}^{\tilde I}$ but a different unmarked restriction on $s'$.
If the restriction is $s[\alpha]=\sigma'_\alpha$,
we rewrite it as the difference of two unmarked restrictions
$s[\alpha]\geq\sigma'_\alpha$ and $s[\alpha]\geq\sigma'_\alpha+1$.
We extend this argument to all $\alpha$ with marked conditions
by using the inclusion-exclusion principle, which will be explained below.

Let
$$
I'=I'_1\sqcup\dots\sqcup I'_h
$$
be the decomposition considered in the proof of Lemma \ref{J-summation}.
Denote by $\tilde I_{\rm right}$ the set of integers ${\rm max}(I'_i)$
for $1\leq i\leq h$ except when it is equal to $k$. This is exactly the
set of $\alpha$ such that the restriction on $s'[\alpha]$ is marked
in $\tilde R_{m',n'}^{\tilde I}$.
For each $K\subset\tilde I_{\rm right}$ we denote by $\tilde I_K$
the subset obtained {}from $\tilde I$ by replacing all the elements
$u\in K$ with $u+1$. Because of (\ref{POSITIVITY}) we have
$$
\sigma'(\tilde I_K)_\alpha=
\begin{cases}
\sigma'(\tilde I)_\alpha+1\quad&\text{if $\alpha\in K\subset[l_1+1,k]$;}\\
\sigma'(\tilde I)_\alpha\quad&\text{otherwise}.
\end{cases}
$$

We set
\begin{eqnarray}
R_{m'.n'}^{\tilde I,\tilde I_K}&=&
\{(\mu',r',\nu',s')\in R_{m',n'};
\nonumber\\
&&r'[\alpha]\geq^{\varepsilon(\tilde I)_\alpha}\rho'(\tilde I)_\alpha
\quad(1\leq\alpha\leq k-1),\nonumber\\
&&s'[\alpha]\geq\sigma'(\tilde I_K)_\alpha\quad(1\leq\alpha\leq k-1)\}.
\end{eqnarray}

For a subset $R\subset R_{m',n'}$
we denote by $1_{{\rm supp}(R)}$ its support function, i.e.,
$1_{{\rm supp}(R)}(x)=1$ if $x\in R$ and $1_{{\rm supp}(R)}(x)=0$
if $x\not\in R$. The inclusion-exclusion principle tells us that
\be
1_{{\rm supp}(\tilde R_{m',n'}^{\tilde I})}
=\sum_{K\subset \tilde I_{\rm right}}(-1)^{\card(K)}
1_{{\rm supp}(R_{m',n'}^{\tilde I,\tilde I_K})}.
\ee
Thus we obtain
\bea\label{right sum}
1_{{\rm supp}(
\bigsqcup_{\tilde I}\tilde R_{m',n'}^{\tilde I}
)}=
\sum_{\tilde I}\sum_{K\subset\tilde I_{\rm right}}(-1)^{\card(K)}
1_{{\rm supp}(R_{m',n'}^{\tilde I,\tilde I_K})}.
\label{right-sum}
\ena

Denote by $\tilde I_{\rm left}$ the set of integers ${\rm min}(I'_i)$ for
$1\leq i\leq h$ except when it is equal to $l_1+1$. This is exactly the set of
$\alpha\in[u_{a+1},k]$ such that  the restriction on
$r'[\alpha-1]$ is marked in  $\overline R_{m',n'}^{\tilde I}$.
For each $K\subset\tilde I_{\rm left}$ we denote by $\tilde I^K$
the subset obtained {}from $\tilde I$ by replacing all the elements
$u\in K$ with $u-1$. Because of (\ref{RHO'2}) we have
$$
\rho'(\tilde I^K)_\alpha=
\begin{cases}
\rho'(\tilde I)_\alpha+1\quad&\text{if $\alpha+1\in K\subset[l_1+1,k]$};\\
\rho'(\tilde I)_\alpha\quad&\text{otherwise}.
\end{cases}
$$

Denote by $\tilde{\cal I}$ the set of $\tilde I$ satisfying (\ref{TIL-I}).
Note that
\begin{eqnarray}
\{(\tilde I,\tilde I_K);\tilde I\in\tilde{\cal I},K\in\tilde I_{\rm right}\}
&=&\{(\tilde I^K,\tilde I);\tilde I\in\tilde{\cal I},K\in\tilde I_{\rm left}\}.
\end{eqnarray}
Therefore, we have
\begin{equation}
(\ref{right-sum})=
\sum_{\tilde I}\sum_{K\subset\tilde I_{\rm left}}(-1)^{\card(K)}
1_{{\rm supp}(R_{m',n'}^{\tilde I^K,\tilde I})}.
\end{equation}

The inclusion-exclusion principle again tells us that
\be
\sum_{K\subset \tilde I_{\rm left}}(-1)^{\card(K)}
1_{{\rm supp}(R_{m',n'}^{\tilde I^K,\tilde I})}
=1_{{\rm supp}(\overline R_{m',n'}^{\tilde I})}
\ee
Therefore, we obtain (\ref{REVERSE}).
\end{proof}

Now we finish the case $l_1+c<k$.

Set $\tilde a=a+k-l_1-c$.
For $I\subset\{1,\dots,k\}$ such that $I=\{u_1,\dots,u_{\tilde a}\}$
($u_1<\dots<u_{\tilde a}$) we define the closure of $I$ by
\begin{equation}
\overline I=
I\backslash\{u_{a+1},\dots,u_{a+d}\}\sqcup[l_1+1,l_1+d]
\end{equation}
where
\begin{equation}
d=
\begin{cases}
0\quad&\text{if $u_{a+1}>l_1+1$};\\
{\rm max}\{i;u_{a+i}\leq l_1+i\}\quad&\text{if $u_{a+1}\leq l_1+1$}.
\end{cases}
\label{DEFD}
\end{equation}
Note that if $d\not=0$ and $\alpha\leq l_1+d$, then
\[
\sigma'(I)_\alpha=\bigl((\kappa[l_1-a+1.k-c]-\kappa(I)\bigr)^+_\alpha=0.
\]
Therefore, we have
\begin{equation}
\sigma'(I)=\sigma'(\overline I).
\label{SIGEQ}
\end{equation}
The set $I$ satisfies the condition $u_{a+1}\geq l_1+1$ if and only if
$\overline I=I$.
\begin{lemma}
Suppose that $\tilde I$ satisfies $(\ref{TIL-I})$. Consider $I$
with its closure equal to $\tilde I$. The subset $R^I_{m',n'}$ is
defined by $(\ref{I})$, in which the definitions
$(\ref{COMP1})$ for $\rho'(I)$ and $(\ref{COMP2})$ for $\sigma'(I)$
are used. We use $l'_1,l'_2,l'_3$ given by $(\ref{THATCASE})$ 
with $a$ replaced  by $\tilde a$. Then, the formulas for $l'_1,l'_2,l'_3$,
$\rho'(I)$ and $\sigma'(I)$ become exactly equal to those used in
$\tilde R^{\tilde I}_{m',n'}$. With this understanding, we have
\begin{equation}
\bigsqcup_{I:\overline I=\tilde I}R^I_{m',n'}=\tilde R^{\tilde I}_{m',n'}.
\label{LASTEQ}
\end{equation}
\label{LASTLEMMA}
\end{lemma}
\begin{proof}
Suppose that $d$ is given by (\ref{DEFD}) for
$\tilde I=\{u_1,\dots,u_{\tilde a}\}$. If $d=0$, then $\overline I=\tilde I$
implies $I=\tilde I$, and (\ref{LASTEQ}) is obvious.
If $d\geq1$, the union in the left hand side is over $I$ such that
\[
I=\tilde I\backslash[l_1+1,l_1+d]\sqcup I'
\]
where
\[
[u_a+1,u_a+d]\geq I'\geq [l_1+1,l_1+d].
\]
Since $\sigma'(I)=\sigma'(\tilde I)$,
by a similar argument as the proof of Lemma \ref{LEM-LOWSUM},
we take the union of $R^I_{m',n'}$ and obtain $\tilde R^{\tilde I}_{m',n'}$.
\end{proof}

\begin{lemma}
We have
\begin{equation}\label{LAST}
\bigsqcup_{\card(\tilde I)=a+k-l_1-c\atop u_{a+1}\geq l_1+1}
\tilde R^{\tilde I}_{m',n'}=R_{m',n'}[l'_1,l'_2,l'_3].
\end{equation}
\end{lemma}
\begin{proof}
This is a consequence of Lemma \ref{N-1SUM} (with $a$ replaced by $\tilde a$)
and Lemma \ref{LASTLEMMA}.
\end{proof}
In conclusion, we have
\begin{proposition}\label{CONCL2}
\begin{equation}
R_{m',n'}[l'_1,l'_2,l'_3]=
\bigsqcup_{\card(I)=a,\card(J)=b}
R_{m',n'}[l_1]^{I,J}.
\label{UPSUM-2}
\end{equation}
\end{proposition}
\begin{proposition}
\begin{equation}
R^{(M,N-1)}_{m',n'}[l'_1,l'_2,l'_3]=
\bigsqcup_{\card(I)=a,\card(J)=b}
R^{(M,N-1)}_{m',n'}[l_1]^{I,J}.
\end{equation}
\end{proposition}

\end{document}